\documentclass[12pt]{amsart}
\usepackage[margin=1in]{geometry}
\usepackage{shor-research}
\usepackage{shortcuts}
\usepackage[numbers]{natbib}

\title{Reflective numerical semigroups}
\author{Caleb M.\ Shor}
\email{cshor@wne.edu}
\address{Western New England University, Springfield, MA, USA}
\begin{document}
\begin{abstract}
    We define a reflective numerical semigroup of genus $g$ as a numerical semigroup that has a certain reflective symmetry when viewed within $\Z$ as an array with $g$ columns. Equivalently, a reflective numerical semigroup has one gap in each residue class modulo $g$. In this paper, we give an explicit description for all reflective numerical semigroups. With this, we can describe the reflective members of well-known families of numerical semigroups as well as obtain formulas for the number of reflective numerical semigroups of a given genus or given Frobenius number.
\end{abstract}
\maketitle

\section{Introduction}
Let $\N_0$ denote the set of non-negative integers, a monoid under addition, and let $\N$ denote the set of positive integers. A numerical semigroup $S$ is a submonoid of $\N_0$ with finite complement. We denote the complement by $\HH(S)=\N_0\setminus S$. Elements of $\HH(S)$ are called gaps of $S$. The genus of $S$, denoted $\genus(S)$, is the number of gaps of $S$. When the genus is positive, the Frobenius number of $S$, denoted $\frob(S)$, is the largest gap of $S$. It happens that every numerical semigroup $S$ can be written as the set of all finite $\N_0$-linear combinations of finitely many natural numbers $a_1,a_2,\dots,a_k$ with $\gcd(a_1,a_2,\dots,a_k)=1$. When this occurs, we write 
\[
    S
    =\left\langle a_1,a_2,\dots,a_k\right\rangle
    =\left\{n_1a_1+n_2a_2+\dots+n_ka_k : n_1,n_2,\dots,n_k\in\N_0\right\}.
\]

In this paper, we are interested in visualizing symmetry in numerical semigroups. It will be helpful to think of a numerical semigroup as a subset of $\Z$. Given a numerical semigroup $S$, we can represent $\Z$ with a 2-dimensional array of boxes and then mark elements of $S$ within it.

As an example, let $S_1=\langle 3, 11\rangle$, a numerical semigroup with genus $\genus(S_1)=10$ and Frobenius number $\frob(S_2)=19$. All integers greater than 19 are in $S_1$, and all integers below 0 are not in $S_1$. To visualize $S_1$, it suffices to view the portion of $S_1$ that lies in the interval $\ints{0,19}$. See Figure~\ref{fig:semigroup-3-11}.

\begin{figure}
    \centering
    \setlength{\unitlength}{.5\linewidth}%
	\mbox{
    \begin{picture}(1, .2)
      \divide\unitlength by 10\relax
      \textcolor{gray}{%
        \put(0, 0){\makerule{1}{1}}%
        \put(3, 0){\makerule{1}{1}}%
        \put(6, 0){\makerule{1}{1}}%
        \put(9, 0){\makerule{1}{1}}%
        \put(1, 1){\makerule{1}{1}}%
        \put(2, 1){\makerule{1}{1}}%
        \put(4, 1){\makerule{1}{1}}%
        \put(5, 1){\makerule{1}{1}}%
        \put(7, 1){\makerule{1}{1}}%
        \put(8, 1){\makerule{1}{1}}%
      }%
      \setcounter{cntA}{-1}%
      \multiput(0, 0)(0, 1){2}{%
        \stepcounter{cntA}%
        \setcounter{cntB}{-1}%
        \multiput(0, 0)(1, 0){10}{%
          \stepcounter{cntB}%
          \makebox(1, 1){\the\numexpr\value{cntA}*10+\value{cntB}\relax}%
        }%
      }   
      \multiput(0, 0)(0, 1){3}{\line(1, 0){10}}
      \multiput(0, 0)(1, 0){11}{\line(0, 1){2}}
    \end{picture}%
    }
    \caption{Visualization of $S_1=\langle 3,11\rangle$, a numerical semigroup of genus $\genus(S_1)=10$, within $\ints{0,19}$. Elements of $S_1$ are shaded. 
    }
    \label{fig:semigroup-3-11}
\end{figure}

In general, all of the gaps of a numerical semigroup $S$ of genus $g=\genus(S)\ge1$ lie in the interval of integers $\ints{0,2g-1}$. (See Proposition~\ref{prop:frob-genus-bounds}.) As a result, we need only think about numerical semigroup elements in the interval $\ints{0,2g-1}$ for our visualizations. We will do so as we have in Figure~\ref{fig:semigroup-3-11}, using boxes in 2 rows and $g$ columns, with boxes for 0 through $g-1$ from left to right in the bottom row, and boxes for $g$ through $2g-1$ from left to right in the top row.

A numerical semigroup $S$ is called \emph{symmetric} when the following occurs: for all $z\in\Z$, exactly one of $z$ and $\frob(S)-z$ is in $S$. Symmetric numerical semigroups have been studied extensively. (See, e.g., \cite{Kunz1970} and \cite[Chapter 3]{RosalesGarciaSanchez09} for connections with one-dimensional local rings.) Our interest here is to visualize their symmetry. For $S_1=\langle 3,11\rangle$, one can verify that $S_1$ is symmetric. For any $S$, let $\rho$ denote a $180^\circ$ rotation about the center of the visualization of $S$. Applying $\rho$ to $S_1$, we see shaded boxes go to unshaded boxes and vice versa. We visualize $S_1$ and its image under $\rho$ in Figure~\ref{fig:semigroup-rotate-1}. For $S_2=\langle 3,13\rangle$, one can verify that $S_2$ is symmetric and that $\rho$ has the same effect, sending shaded boxes to unshaded boxes and vice versa. We visualize $S_2$ and its image under $\rho$ in Figure~\ref{fig:semigroup-rotate-2}. In general, with the way we are visualizing a numerical semigroup $S$, $\rho$ amounts to swapping the boxes of $z$ and $\frob(S)-z$ for all $z\in\ints{0,2g-1}$. Note that if we drew all of $\Z$ instead of just $\ints{0,2g-1}$, this rotation would swap boxes above $2g-1$, which are in $S$, with boxes below $0$, which are not in $S$.

\begin{figure}[t]
    \centering
    \begin{minipage}[c]{.35\textwidth}
    \mbox{}\hfill
    \setlength{\unitlength}{.8\linewidth}%
	\mbox{
    \begin{picture}(1, .2)
      \divide\unitlength by 10\relax
      \textcolor{gray}{%
        \put(0, 0){\makerule{1}{1}}%
        \put(3, 0){\makerule{1}{1}}%
        \put(6, 0){\makerule{1}{1}}%
        \put(9, 0){\makerule{1}{1}}%
        \put(1, 1){\makerule{1}{1}}%
        \put(2, 1){\makerule{1}{1}}%
        \put(4, 1){\makerule{1}{1}}%
        \put(5, 1){\makerule{1}{1}}%
        \put(7, 1){\makerule{1}{1}}%
        \put(8, 1){\makerule{1}{1}}%
      }%
      \multiput(0, 0)(0, 1){3}{\line(1, 0){10}}
      \multiput(0, 0)(1, 0){11}{\line(0, 1){2}}
    \end{picture}%
    }\hfill\mbox{}
    \end{minipage}
    \ \  $\overset{\rho}{\longrightarrow}$
    \begin{minipage}[c]{.35\textwidth}
    \mbox{}\hfill
    \setlength{\unitlength}{.8\linewidth}%
	\mbox{
    \begin{picture}(1, .2)
      \divide\unitlength by 10\relax
      \textcolor{gray}{%
        \put(1, 0){\makerule{1}{1}}%
        \put(2, 0){\makerule{1}{1}}%
        \put(4, 0){\makerule{1}{1}}%
        \put(5, 0){\makerule{1}{1}}%
        \put(7, 0){\makerule{1}{1}}%
        \put(8, 0){\makerule{1}{1}}%
        \put(0, 1){\makerule{1}{1}}%
        \put(3, 1){\makerule{1}{1}}%
        \put(6, 1){\makerule{1}{1}}%
        \put(9, 1){\makerule{1}{1}}%
      }%
      \multiput(0, 0)(0, 1){3}{\line(1, 0){10}}
      \multiput(0, 0)(1, 0){11}{\line(0, 1){2}}
    \end{picture}%
    }\hfill\mbox{}
    \end{minipage}
    \caption{Visualization of $S_1=\langle 3,11\rangle$, a numerical semigroup of genus $\genus(S_1)=10$, within $\ints{0,19}$ along with its image under $\rho$, a $180^\circ$ rotation.
    }
    \label{fig:semigroup-rotate-1}
\end{figure}
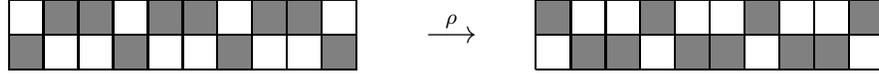
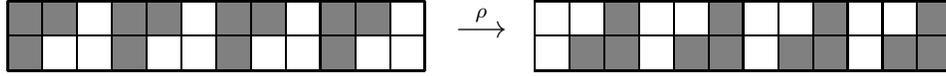
\begin{figure}[t]
    \centering
    \begin{minipage}[c]{.35\textwidth}
    \mbox{}\hfill
    \setlength{\unitlength}{.96\linewidth}%
	\mbox{
    \begin{picture}(1, .2)
      \divide\unitlength by 12\relax
      \textcolor{gray}{%
        \put(0, 0){\makerule{1}{1}}%
        \put(3, 0){\makerule{1}{1}}%
        \put(6, 0){\makerule{1}{1}}%
        \put(9, 0){\makerule{1}{1}}%
        \put(0, 1){\makerule{1}{1}}%
        \put(1, 1){\makerule{1}{1}}%
        \put(3, 1){\makerule{1}{1}}%
        \put(4, 1){\makerule{1}{1}}%
        \put(6, 1){\makerule{1}{1}}%
        \put(7, 1){\makerule{1}{1}}%
        \put(9, 1){\makerule{1}{1}}%
        \put(10, 1){\makerule{1}{1}}%
      }%
      \multiput(0, 0)(0, 1){3}{\line(1, 0){12}}
      \multiput(0, 0)(1, 0){13}{\line(0, 1){2}}
    \end{picture}%
    }\end{minipage}
    \ \  $\overset{\rho}{\longrightarrow}$
    \begin{minipage}[c]{.35\textwidth}
    \mbox{}\hfill
    \setlength{\unitlength}{.96\linewidth}%
	\mbox{
    \begin{picture}(1, .2)
      \divide\unitlength by 12\relax
      \textcolor{gray}{%
        \put(1, 0){\makerule{1}{1}}%
        \put(2, 0){\makerule{1}{1}}%
        \put(4, 0){\makerule{1}{1}}%
        \put(5, 0){\makerule{1}{1}}%
        \put(7, 0){\makerule{1}{1}}%
        \put(8, 0){\makerule{1}{1}}%
        \put(10, 0){\makerule{1}{1}}%
        \put(11, 0){\makerule{1}{1}}%
        \put(2, 1){\makerule{1}{1}}%
        \put(5, 1){\makerule{1}{1}}%
        \put(8, 1){\makerule{1}{1}}%
        \put(11, 1){\makerule{1}{1}}%
      }%
      \multiput(0, 0)(0, 1){3}{\line(1, 0){12}}
      \multiput(0, 0)(1, 0){13}{\line(0, 1){2}}
    \end{picture}%
    }\end{minipage}
    \caption{Visualization of $S_2=\langle 3,13\rangle$, a numerical semigroup of genus $\genus(S_2)=12$, within $\ints{0,23}$ along with its image under $\rho$, a $180^\circ$ rotation.
    }
    \label{fig:semigroup-rotate-2}
\end{figure}

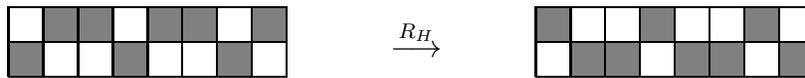
\begin{figure}[b]
    \centering
    \begin{minipage}[c]{.35\textwidth}
    \mbox{}\hfill
    \setlength{\unitlength}{.64\linewidth}%
	\mbox{
    \begin{picture}(1, .2)
      \divide\unitlength by 8\relax
      \textcolor{gray}{%
        \put(0, 0){\makerule{1}{1}}%
        \put(3, 0){\makerule{1}{1}}%
        \put(6, 0){\makerule{1}{1}}%
        \put(1, 1){\makerule{1}{1}}%
        \put(2, 1){\makerule{1}{1}}%
        \put(4, 1){\makerule{1}{1}}%
        \put(5, 1){\makerule{1}{1}}%
        \put(7, 1){\makerule{1}{1}}%
      }%
      \multiput(0, 0)(0, 1){3}{\line(1, 0){8}}
      \multiput(0, 0)(1, 0){9}{\line(0, 1){2}}
    \end{picture}%
    }\hfill\mbox{}
    \end{minipage}
    \ \  $\overset{\reflect}{\longrightarrow}$
    \begin{minipage}[c]{.35\textwidth}
    \begin{center}
    \mbox{}\hfill
    \setlength{\unitlength}{.64\linewidth}%
    \mbox{
    \begin{picture}(1, .2)
      \divide\unitlength by 8\relax
      \textcolor{gray}{%
        \put(1, 0){\makerule{1}{1}}%
        \put(2, 0){\makerule{1}{1}}%
        \put(4, 0){\makerule{1}{1}}%
        \put(5, 0){\makerule{1}{1}}%
        \put(7, 0){\makerule{1}{1}}%
        \put(0, 1){\makerule{1}{1}}%
        \put(3, 1){\makerule{1}{1}}%
        \put(6, 1){\makerule{1}{1}}%
      }%
      \multiput(0, 0)(0, 1){3}{\line(1, 0){8}}
      \multiput(0, 0)(1, 0){9}{\line(0, 1){2}}
    \end{picture}%
    }\hfill\mbox{}
    \end{center}
    \end{minipage}
    
    \caption{Visualization of $S_3=\langle 3,10,17\rangle$, a numerical semigroup of genus $\genus(S_3)=8$, within $\ints{0,15}$ along with its image under $\reflect$, a horizontal reflection.
    }
    \label{fig:semigroup-reflect-1}
\end{figure}

\begin{figure}[b]
    \centering
    \begin{minipage}[c]{.35\textwidth}
    \mbox{}\hfill
    \setlength{\unitlength}{.56\linewidth}%
	\mbox{
    \begin{picture}(1, .2)
      \divide\unitlength by 7\relax
      \textcolor{gray}{%
        \put(0, 0){\makerule{1}{1}}%
        \put(4, 0){\makerule{1}{1}}%
        \put(1, 1){\makerule{1}{1}}%
        \put(2, 1){\makerule{1}{1}}%
        \put(3, 1){\makerule{1}{1}}%
        \put(5, 1){\makerule{1}{1}}%
        \put(6, 1){\makerule{1}{1}}%
      }%
      \multiput(0, 0)(0, 1){3}{\line(1, 0){7}}
      \multiput(0, 0)(1, 0){8}{\line(0, 1){2}}
    \end{picture}%
    }\hfill\mbox{}
    \end{minipage}
    \ \  $\overset{\reflect}{\longrightarrow}$
    \begin{minipage}[c]{.35\textwidth}
    \begin{center}
    \mbox{}\hfill
    \setlength{\unitlength}{.56\linewidth}%
    \mbox{
    \begin{picture}(1, .2)
      \divide\unitlength by 7\relax
      \textcolor{gray}{%
        \put(1, 0){\makerule{1}{1}}%
        \put(2, 0){\makerule{1}{1}}%
        \put(3, 0){\makerule{1}{1}}%
        \put(5, 0){\makerule{1}{1}}%
        \put(6, 0){\makerule{1}{1}}%
        \put(0, 1){\makerule{1}{1}}%
        \put(4, 1){\makerule{1}{1}}%
      }%
      \multiput(0, 0)(0, 1){3}{\line(1, 0){7}}
      \multiput(0, 0)(1, 0){8}{\line(0, 1){2}}
    \end{picture}%
    }\hfill\mbox{}
    \end{center}
    \end{minipage}
    
    \caption{Visualization of $S_4=\langle 4,9,10,15\rangle$, a numerical semigroup of genus $\genus(S_4)=7$, within $\ints{0,13}$ along with its image under $\reflect$, a horizontal reflection.
    }
    \label{fig:semigroup-reflect-2}
\end{figure}
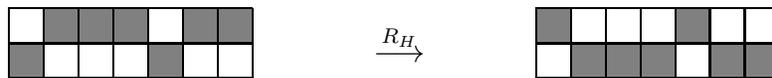

The visualization of $S_1$ has a further symmetry that is not present for $S_2$: reflection across the horizontal line that separates the two rows, which we denote by $\reflect$. When a numerical semigroup has this symmetry, we say the numerical semigroup is \emph{reflective}. Thus, $S_1$ is reflective and $S_2$ is not reflective. In Figure~\ref{fig:semigroup-reflect-1} and Figure~\ref{fig:semigroup-reflect-2}, we visualize the numerical semigroups $S_3=\langle 3,10,17\rangle$ and $S_4=\langle 4,9,10,15\rangle$ and their images under $R_H$. Both $S_3$ and $S_4$ are reflective. Neither $S_3$ nor $S_4$ is symmetric, though $S_3$ is pseudo-symmetric (a condition which is defined in Section~\ref{sec:ME-props}).

As we will see (Prop.~\ref{prop:reflective-equiv-defn}), a  numerical semigroup $S$ is reflective precisely when it has one gap in each residue class modulo $\genus(S)$. In the language of \cite{Shor22}, we say such a numerical semigroup is equidistributed modulo $\genus(S)$.

The purpose of this paper is to describe all reflective numerical semigroups. We will give an explicit description (Proposition~\ref{prop:ED iff cases}) of any reflective numerical semigroup in terms of its genus and multiplicity. (The multiplicity of a numerical semigroup is defined in Section~\ref{sec:prelims}.) Among numerical semigroups of positive genus, there is one reflective numerical semigroup for each combination of genus $g\ge1$ and multiplicity $a$ with $2\le a\le g+1$ and $a\nmid g$, and there are no others. We will identify the particular members of a few frequently-studied families of numerical semigroups (such as those with embedding dimension 2, those generated by generalized arithmetic sequences, and free numerical semigroups) that have this reflective property.

Finally, we will describe all, and count the number of, reflective numerical semigroups of a particular genus or Frobenius number. Let  $\tau(n)$ (resp., $\tau_e(n)$) denote the number of positive divisors (resp., positive even divisors) of a positive integer $n$. Then the number of reflective numerical semigroups of genus $g$ is $g+1-\tau(g)$ (by Prop.~\ref{prop:number-of-ME-semigroups-by-genus}), and the number of reflective numerical semigroups of Frobenius number $F$ is
\[
    1-\tau_e(F)+\sum\limits_{k=2}^F(-1)^k\Floor{\frac{F}{k}} 
    = (1-\log2)F+\bigO(\sqrt{F})
\]
(by Cor.~\ref{cor:nF-exact} and Prop.~\ref{prop:nF-asymp}). We also compute the number of numerical semigroups of genus $g$ that are both reflective and symmetric (Prop.~\ref{prop:num-symm-MEs}).

\subsection{Organization}
Our work is organized as follows. In Section \ref{sec:prelims}, we provide the necessary background needed for this paper. In Section \ref{sec:ED mod g}, we find an explicit description for a reflective numerical semigroup of genus $g$ and multiplicity $a$. We then compute some important properties, such as the Frobenius number, the embedding dimension, the minimal generating set, and an Ap\'ery set of such a numerical semigroup. In Section \ref{sec:ME-props}, we describe a few well-known families of numerical semigroups and identify the members of those families that are reflective. We conclude this section by showing that Wilf's conjecture holds for reflective numerical semigroups. In Section \ref{sec:counting-ME-semigroups}, we count the number of reflective numerical semigroups of a particular genus or Frobenius number, concluding with the asymptotic result mentioned above.

\section{Preliminaries}\label{sec:prelims}
In this paper, we will work with intervals of integers. For $a,b\in\Z$ with $a\le b$, we will therefore let $\ints{a,b}$ denote the interval of integers from $a$ to $b$. (That is, $\ints{a,b}=\{n\in\Z : a\le n\le b\}$.)  
The cardinality of a finite set $A$ is denoted $\#A$. For any real number $x$, we have the floor of $x$ and the ceiling of $x$ defined, respectively, in the usual way by $\left\lfloor x\right\rfloor =\max\{z\in\Z : z\le x\}$ and $\left\lceil x\right\rceil = \min\{z\in\Z : z\ge x\}$. For $a,b\in\Z$, we write $a\mid b$ to mean $b=ak$ for some $k\in\Z$, and similarly $a\nmid b$ to mean there is no such $k\in\Z$ such that $b=ak$.

\subsection{Numerical semigroups}
In the introduction, we defined a numerical semigroup $S$ along with its set of gaps $\HH(S)$, its genus $\genus(S)$, and its Frobenius number $\frob(S)$. In this section, we will give more preliminary information about numerical semigroups. For an in-depth treatment -- including motivation, use, and proofs -- see \cite{RosalesGarciaSanchez09}.

We begin with relations between the genus and Frobenius number of $S$. Since gaps are positive integers and $\frob(S)$ is the largest gap, we have $\HH(S)\subseteq\ints{1,\frob(S)}$. Thus, $\genus(S)\le\frob(S)$, a lower bound for $\frob(S)$.

For an upper bound, note that if the sum of two non-negative integers is equal to $\frob(S)$, then at least one of the two integers must be a gap. (Otherwise, $\frob(S)$ is a sum of elements of $S$, implying $\frob(S)\in S$, a contradiction.) Hence, at least half of the integers in $\ints{0,\frob(S)}$ are gaps. This means $\genus(S)\ge (\frob(S)+1)/2$.

Combining the results of the two previous paragraphs, we obtain the following.

\begin{prop}
\label{prop:frob-genus-bounds}
    For any numerical semigroup $S$, 
    \[
        \genus(S)\le \frob(S)\le 2\genus(S)-1
    \]
    and $\HH(S)\subseteq\ints{1,2\genus(S)-1}$.
\end{prop}

We now describe generating sets of numerical semigroups. For any set $A\subseteq\N_0$, let $\langle A\rangle$ denote the set of all finite $\N_0$-linear combinations of elements of $A$. That is, let 
\[
    \langle A\rangle 
    = \left\{\sum\limits_{i=1}^k n_ia_i : k\in\N_0, n_i\in\N_0, a_i\in A \right\},
\] 
the $\N_0$-module generated by $A$. It is a standard result that $\langle A\rangle$ is a numerical semigroup if and only if $\gcd(A)=1$. When this occurs, we say $A$ is a \emph{generating set} for $\langle A\rangle$. As was mentioned in the introduction, every numerical semigroup $S$ has a generating set. If $S=\langle A\rangle$, then the smallest positive element of $S$ is the smallest positive element of $A$, and we call this number the \emph{multiplicity} of $S$, denoted $\mult(S)$.

A generating set $A$ is \emph{minimal} if $\langle A\setminus\{a\}\rangle\ne\langle A\rangle$ for all $a\in A$. Equivalently, this occurs if no element of $A$ is a $\N_0$-linear combination of the other elements of $A$. For any numerical semigroup $S$, there is a unique minimal (with respect to set inclusion) generating set $A$ for which $S=\langle A\rangle$. Such a minimal generating set is necessarily finite. The \emph{embedding dimension} of $S$, denoted $\embdim(S)$, is the cardinality of this minimal generating set.

An important tool for use with numerical semigroups is the Ap\'ery set of a numerical semigroup relative to a nonzero element.

\begin{defn}[\cite{Apery1946}]
    For $S$ a numerical semigroup and $0\ne t\in S$, the \emph{Ap\'ery set of $S$ relative to $t$} is \[\Ap(S;t)=\left\{s\in S : s-t\not\in S\right\}.\] 
\end{defn}

The set $\Ap(S;t)$ consists of smallest element of $S$ in each congruence class modulo $t$. It follows that $\#\Ap(S;t)=t$. Viewing $\Ap(S;t)$ in this way, we see that $\frob(S)$, which is the largest gap of $S$, is simply the largest element of $\Ap(S;t)$ minus $t$. One can compute $\genus(S)$ by counting the number of gaps in each congruence class modulo $t$ and simplifying the sum. We state these results in the following.

\begin{prop}[{\cite[Proposition 2.12]{RosalesGarciaSanchez09}}]
    For any nonzero $t\in S$, 
    \[
        \frob(S)
        =-t+\max\left(\Ap(S;t)\right)
        \text{ and } 
        \genus(S)
        =\frac{1-t}{2}+\frac{1}{t}\sum\limits_{w\in\Ap(S;t)}w.
    \]
\end{prop}

Furthermore, if two numerical semigroups have the same Ap\'ery set relative to a particular element, then the two semigroups are necessarily equal.
\begin{prop}[{\cite[Lemma 5.4]{Shor22}}]\label{prop:same-apery-same-semigroup}
    Let $S,T$ be numerical semigroups. For nonzero $t\in S\cap T$, we have $\Ap(S;t)=\Ap(T;t)$ if and only if $S=T$.
\end{prop}

\subsection{Reflective numerical semigroups and equidistributed sets}
In the introduction, we defined a reflective numerical semigroup $S$ in terms of a certain visual symmetry. For $g=\genus(S)$, the gaps of $S$ are contained in the interval $\ints{0,2g-1}$. To visualize $S$, we can focus on elements of $S$ in the interval $\ints{0,2g-1}$. We make a rectangle out of boxes in two rows: the bottom row $\ints{0,g-1}$ and the top row $\ints{g,2g-1}$. We shade elements of $S$. A numerical semigroup is reflective if, when we reflect across the horizontal line separating the rows, semigroup elements go to semigroup gaps and vice versa.  This reflection is denoted $\reflect$.

\begin{defn}[Reflective numerical semigroup]\label{defn:reflective-v1}
    A numerical semigroup $S$ of genus $g=\genus(S)$ is \emph{reflective} if, for all $z\in\ints{0,g-1}$, exactly one of $z$ and $z+g$ is in $S$.
\end{defn}

\begin{remark}
    If $\genus(S)=0$, then $S=\N_0$ and by the above definition, $S$ is vacuously reflective. For the rest of this paper, we will suppose $g\ge1$.
\end{remark}

The standard definition of a symmetric numerical semigroup (which we stated in the introduction) holds for all $z\in\Z$. It would be nice to have something similar for a reflective numerical semigroup. In Figure~\ref{fig:semigroup-4-9-10-15-again}, we have visualized $S_4=\langle 4,9,10,15\rangle$ again, this time within a larger subinterval of $\Z$. Thinking about our visualization within $\Z$, all boxes above the rectangle are shaded, and all boxes below are unshaded. Consider the column containing 2 and 9. When we apply our reflection $\reflect$, the following pairs of boxes in this column are mapped to each other: $2$ and $9$; $-5$ and $16$; $-12$ and $23$; etc. Note that the sum of each reflected pair in this column is 11, a constant. If we start with another column and add reflected pairs together, we get a different constant.

For any integer $z$, we can divide $z$ by $g$ to get integers $q$ and $r$ such that $z=qg+r$ with $0\le r<g$. In a reflective numerical semigroup, for each $r$ with $0\le r<g$, exactly one of $r$ and $g+r$ is in $S$. Their sum is $g+2r$. Extending upward and downward by $k$ rows in our picture, we see that exactly one of $-kg+r$ and $(k+1)g+r$ is in $S$ for any $k\in\Z$. When this occurs, $z$ and its reflection sum to $g+2r$. Thus, the reflection of $z$ is $g+2r-z$. We can now extend our definition to all $z\in\Z$, keeping in mind that $r=z-qg=z-\Floor{z/g}g$.

\begin{defn}[Reflective numerical semigroup, extended version]\label{defn:reflective-v2}
    Let $S$ be a numerical semigroup of genus $g=\genus(S)\ge1$. We say $S$ is  \emph{reflective} if, for all $z\in\Z$, exactly one of $z$ and $z+\left(1-2\left\lfloor\dfrac{z}{g}\right\rfloor\right)g$ is in $S$.
\end{defn}
\begin{figure}
\centering
    \setlength{\unitlength}{.5\linewidth}%
	\mbox{
    \begin{picture}(.7, .6)
      \divide\unitlength by 10\relax
      \textcolor{gray}{%
        \put(1, 2){\makerule{1}{1}}%
        \put(5, 2){\makerule{1}{1}}%
        \put(2, 3){\makerule{1}{1}}%
        \put(3, 3){\makerule{1}{1}}%
        \put(4, 3){\makerule{1}{1}}%
        \put(6, 3){\makerule{1}{1}}%
        \put(7, 3){\makerule{1}{1}}%
        \multiput(1,4)(1,0){7}{\makerule{1}{1}}%
        \multiput(1,5)(1,0){7}{\makerule{1}{1}}%
      }%
      \setcounter{cntA}{-1}%
      \multiput(1, 0)(0, 1){6}{%
        \stepcounter{cntA}%
        \setcounter{cntB}{-1}%
        \multiput(0, 0)(1, 0){7}{%
          \stepcounter{cntB}%
          \makebox(1, 1){$\the\numexpr\value{cntA}*7-14+\value{cntB}$\relax}%
        }%
      }   
      \multiput(1, 0)(0, 1){7}{\line(1, 0){7}}
      \multiput(1, 0)(1, 0){8}{\line(0, 1){6}}
      \thinlines 
      \put(0.0,3){\line(1,0){0.1}}
      \put(0.2,3){\line(1,0){0.1}}
      \put(0.4,3){\line(1,0){0.1}}
      \put(0.6,3){\line(1,0){0.1}}
      \put(0.8,3){\line(1,0){0.1}}
      \put(8.0,3){\line(1,0){0.1}}
      \put(8.2,3){\line(1,0){0.1}}
      \put(8.4,3){\line(1,0){0.1}}
      \put(8.6,3){\line(1,0){0.1}}
      \put(8.8,3){\line(1,0){0.1}}
      \put(9.1,3){\vector(0,1){0.5}}
      \put(9.1,3){\vector(0,-1){0.5}}
      \put(9.4,2.9){$\reflect$}
    \end{picture}%
    }
    \caption{Another visualization of $S_4=\langle 4,9,10,15\rangle$, a numerical semigroup of genus $\genus(S_4)=7$, this time within $\ints{-14,27}$. Elements of $S_4$ are shaded. The symmetry $\reflect$ reflects across the (dashed) horizontal line between the rows $\ints{0,6}$ and $\ints{7,13}$. 
    }
    \label{fig:semigroup-4-9-10-15-again}
\end{figure}

We will conclude this section with alternative characterization of reflective numerical semigroups. Let $A\subset\Z$ be a finite set, and let $m\in\mathbb{N}$. We take the definition of equidistribution modulo $m$ and two preliminary results about it from \cite{Shor22}.

\begin{defn}\label{def:n_{r,m}}
    For $r\in\Z$,  let $A_{r,m}=\{a\in A : a\equiv r\pmod*{m}\}$. We define
    \[
        \nn_{r,m}(A) = \#A_{r,m},
    \] 
    the number of elements of $A$ that are congruent to $r$ modulo $m$. The set $A$ is \emph{equidistributed modulo $m$} if $\nn_{r_1,m}(A)=\nn_{r_2,m}(A)$ for all $r_1,r_2\in\ints{1,m}$.
\end{defn}

Put another way, $A$ is equidistributed modulo $m$ if $A$ contains the same number of elements in each congruence class modulo $m$. We now state a few results about sets which are equidistributed modulo $m$.

\begin{lem}[{\cite[Lemma 2.8]{Shor22}}]\label{lem:equidist-candidates}
    Suppose $A$ is equidistributed modulo $m$. Then $m\mid\#A$, and $\#A=m\cdot\nn_{r,m}(A)$ for all $r\in\ints{1,m}$.
\end{lem}

In other words, $A$ can only possibly be equidistributed modulo divisors of $\#A$. We'll say a set $A$ is \emph{maximally equidistributed} if $A$ is equidistributed modulo all positive divisors of $\#A$.

\begin{lem}[{\cite[Proposition 2.10]{Shor22}}]\label{lem:equidist-divisors}
     If $A$ is equidistributed modulo $m$ then $A$ is equidistributed modulo $d$ for each positive divisor $d$ of $m$.
\end{lem}

It follows that a set $A$ is maximally equidistributed if and only if $A$ is equidistributed modulo $\#A$. We may now characterize reflective numerical semigroups as those which have a set of gaps that is equidistributed modulo $m$ for as many $m$ as possible.

\begin{prop}\label{prop:reflective-equiv-defn}
    A numerical semigroup $S$ is reflective if and only if its set of gaps $\HH(S)$ is maximally equidistributed, which occurs if and only if $\HH(S)$ is equidistributed modulo $\genus(S)$.
\end{prop}
\begin{proof}
    By Definition~\ref{defn:reflective-v1}, $S$ is reflective if and only if $\HH(S)$ has one gap in each congruence class modulo $\genus(S)$. Hence, $S$ is reflective if and only if $\HH(S)$ is equidistributed modulo $\genus(S)$, which occurs if and only if $\HH(S)$ is equidistributed modulo all divisors of $\genus(S)$.
\end{proof}

\section{Generating elements of reflective numerical semigroups}\label{sec:ED mod g}

The material in this section is split into two parts. In Section~\ref{subsec:explicit-ME-description}, we will explicitly describe all reflective numerical semigroups. In Section~\ref{subsec:ME-props}, we will determine the Frobenius number, the embedding dimension, and the Ap\'ery set relative to the smallest generating element of any reflective numerical semigroup.

\subsection{General form of reflective numerical semigroups}\label{subsec:explicit-ME-description}

We will now work toward an explicit description of each reflective numerical semigroup $S$. We start with a lemma.

\begin{lem}\label{lem:n in S iff g=0}
    Suppose the set of gaps of a numerical semigroup $S$ is equidistributed modulo $m$ for some $m\in\mathbb{N}$. Then $m\in S$ if and only if $\genus(S)=0$.
\end{lem}
\begin{proof}
    Let $m\in\N$ and suppose the set of gaps of $S$ is equidistributed modulo $m$. If $m\in S$, then all integral multiples of $m$ are in $S$. Therefore there are no gaps that are 0 modulo $m$. Since $S$ is equidistributed modulo $m$, there can be no gaps in any other congruence class modulo $m$. Thus, $\genus(S)=0$.

    For the converse, if $\genus(S)=0$, then $S=\mathbb{N}_0$. Thus, $m\in S$.
\end{proof}
\begin{remark}\label{rmk:g in s iff g=0}
    If $S$ is reflective, then $\genus(S)\in S$ if and only if $\genus(S)=0$. Our main takeaway is that, if $\genus(S)>0$, then since $\mult(S)\in S$ and $S$ is closed under addition, we must have $\mult(S)\nmid\genus(S)$.
\end{remark}

We are almost ready to describe a generating set for any reflective numerical semigroup (Proposition~\ref{prop:ED iff cases}). The following lemma will be useful for the proof of our main result, allowing us to assume that elements of a certain numerical semigroup can be written in a particular form.

\begin{lem}\label{lem:toward-minimality}
    For $g\in\N$ and $a\in\ints{2,g+1}$ with $a\nmid g$, consider the set 
    \[A=\{a\}\cup\ints{g+1,g+(a-1)}.\]
    If $z\in\langle A\rangle$ and $g<z\le 2g+1$, then $z=la+(g+i)$ for some $l\in\N_0$ and some $i\in\ints{1,a-1}$.
\end{lem}
\begin{proof}
    Suppose $z\in\langle A\rangle$ and $g<z\le 2g+1$. Then 
    \[
        z
        =la+\sum\limits_{i=1}^{a-1}c_i(g+i)
    \] 
    for some $l,c_1,\dots,c_{a-1}\in\N_0$. If $\sum c_i\ge2$, then $z\ge2(g+1)=2g+2$, a contradiction. Hence, $\sum c_i$ is 0 or 1. If the sum is 1, then we are done.
    
    If the sum is 0, then $z=da$ for some $d\in\N_0$. By the division algorithm in $\Z$, there are $q,r\in\Z$ with $g=qa+r$ and $0\le r<m$. Since $a\nmid g$, we have $r\ne0$. Then the smallest multiple of $a$ that is greater than $g=qa+r$ is $g+(a-r)=(q+1)a$. Therefore, since $a\mid z$ and $z\ge g$, $z\ge (q+1)a$ and hence $d\ge q+1$. It follows that 
    \[
        z=d a =(d-(q+1))a + (q+1)a,
    \]
    so $z=la+(g+i)$ for $l=d-(q+1)\ge0$ and $i=a-r\in\ints{1,a-1}$, as desired.
\end{proof}

We now present a generating set for any reflective numerical semigroup. Recall that the \emph{multiplicity} of a numerical semigroup $S$, denoted $\mult(S)$, is the minimal nonzero element of $S$.

\begin{prop}\label{prop:ED iff cases}
    Suppose $S$ is a numerical semigroup with $\genus(S)=g\ge1$ and $\mult(S)=a$. Let $q=\left\lfloor g/a\right\rfloor$ and $r=g-aq$. Then $S$ is reflective if and only if $S=\langle A\rangle$ for 
    \[
        A=\{a, 2g+a-r\}\cup \ints{g+1,g+(a-1)},
    \] 
    with $a\in\ints{2,g+1}$ and $a\nmid g$.
\end{prop}
\begin{proof}
    First, note that in defining $q$ and $r$ as above, we can equivalently define them as integers resulting from the division algorithm on $g$ and $a$, for which $g=aq+r$ and $0\le r<a$. Next, since $g\ge1$, $1\not\in S$, so we must have $a\ge2$.

    ($\implies$) 
    Suppose $S$ is reflective. If $r=0$, then $a\mid g$. Since $a\in S$, we have $g\in S$. By Remark \ref{rmk:g in s iff g=0}, this implies $g=0$, contradicting our assumption that $g\ge1$.

    For the remainder of this proof, we can assume $r\ne 0$ (or, equivalently that $a\nmid g$). Let $I_1=\ints{0,g-1}$ and $I_2=\ints{g,2g-1}$. Since $S$ is reflective, for each $z\in I_1$ we have that exactly one of $z$ and $z+g$ is in $S$. Similarly, for each $z\in I_2$, exactly one of $z$ and $z-g$ is in $S$.
    
    We have $\mult(S)=a$. Thus, $1,2,\dots,a-1\not\in S$, implying $\ints{1,a-1}\subset\HH(S)$. There are at least $a-1$ gaps, which implies $g\ge a-1$ and thus $a\in\ints{2,g+1}$. We also have $\{0,a,2a,\dots,qa\}\in S$ because $a\in S$ and $S$ is closed under addition. Since $\ints{1,a-1}\cup\{0,a,2a,\dots,qa\}\subset I_1$, we add $g$ to each and find that $\ints{g+1,g+(a-1)}\subset S$ and $\{g,g+a,g+2a,\dots,g+qa\}\subset\HH(S)$. Since $a\in S$, we can add $a$ repeatedly to elements of $\ints{g+1,g+(a-1)}$ to conclude that all elements of $I_2$ which are not congruent to $g$ modulo $a$ are therefore in $S$. Thus $\HH(S)\cap I_2 = \{g,g+a,g+2a,\dots,g+qa\}$, from which it follows that $S\cap I_1=\{0,a,2a,\dots,qa\}$. Note that the largest gap is $g+qa=2g-r$.
    
    We can therefore generate $S$ with the union of three sets: $\{a\}$, $\ints{g+1,g+(a-1)}$, and $\{2g+a-r\}$. The set $\{a\}$ generates all elements of $S\cap I_1$. The sets $\{a\}$ and $\ints{g+1,g+(a-1)}$ together generate all elements of $S\cap I_2$. Since all integers greater than $2g-1$ are in $S$, and all of the gaps in $I_2$ are congruent to $g$ modulo $a$, we add $a$ to the largest gap to get $g+qa+a=2g+a-r$, our final generating element. This ensures that all integers greater than $2g-1$ are in $S$. Thus, we have $S=\langle A\rangle$ for 
    \[
        A=\{a\}\cup\ints{g+1,g+(a-1)}\cup\{2g+a-r\}.
    \] 
    As mentioned above, $a\in\ints{2,g+1}$ and $a\nmid g$, as desired.
    
    ($\impliedby$) 
    Now, suppose $a\in\ints{2,g+1}$, $a \nmid g$, and $S=\langle A\rangle$ for 
    \[
        A=\{a,2g+a-r\}\cup\ints{g+1,g+(a-1)}.
    \] 
    We wish to show that $\genus(S)=g$ and that $\HH(S)$ is equidistributed modulo $g$. We begin by showing that all of the gaps lie in the interval $\ints{0,2g-1}.$
    
    Let $z\ge 2g$. We wish to show that $z\in S$. First, suppose $z\not\equiv g\pmod*{a}$, so $z\equiv g+i\pmod*{a}$ for some $i\in\ints{1,a-1}\subseteq\ints{1,g}$. Since $z\ge 2g\ge g+i$, $z=(g+i)+la$ for some $l\in\mathbb{Z}$ with $l\ge0$. Thus $z\in S$. Next, suppose $z\equiv g\pmod*{a}$, so $z\equiv 2g+a-r\pmod*{a}$. Thus $z=2g+a-r+la$ for some $l\in\mathbb{Z}$. Solving for $l$, we find that 
    \[
        l=(z-2g-a+r)/a\ge (2g-2g-a+r)/a=(r-a)/a>-1.
    \] 
    Since $l\in\mathbb{Z}$, $l\ge0$ and therefore $z\in S$.
    
    Now that we know $\HH(S)\subset[0,2g-1]$, we again split this interval into $I_1=\ints{0,g-1}$ and $I_2=\ints{g,2g-1}$. We will show that for each $i\in I_1$, exactly one of $i$ and $i+g$ is in $S$. 
    
    We first look at $I_1$. Since $qa<g+i$ for all $i\in\ints{1,a-1}$ and $qa<2g+a-r$, the only generating element of $S$ in $I_1$ is $a$. Thus, $S\cap I_1=\{0,a,\dots,qa\}$, from which it follows that $\HH(S)\cap I_1=\{0\le z\le g-1 : z\not\equiv 0\pmod*{a}\}$.
    
    Now we look at $I_2$. We have $\ints{g+1,g+(a-1)}\subset I_2\cap S$. Let $z\in I_2$. If $z\equiv g+i\pmod*{a}$ for some $i\in\ints{1,a-1}$, then since $z\ge g+i$, $z=(g+i)+la$ for some $l\in\N_0$, which implies $z\in S$. This covers all congruence classes for $z$ except for $g$ modulo $a$. If $z\equiv g\pmod*{a}$, for the sake of contradiction assume that $z\in I_2\cap S$. By Lemma \ref{lem:toward-minimality},  $z=la+(g+i)$ for some $l\in\N_0$ and $i\in\ints{1,a-1}$. However, $z=la+(g+i)\equiv g+i\not\equiv g\pmod*{a}$, which is a contradiction. Thus, if $z\equiv g\pmod*{a}$, then $z\not\in I_2\cap S$. We have therefore found that $S\cap I_2=\{g\le z\le 2g-1 : z\not\equiv g\pmod*{a}\}$ and $\HH(S)\cap I_2=\{g,g+a,g+2a,\dots,g+qa\}$.
    
    Combining the previous two paragraphs, we see that exactly one of $i$, $i+g$ is in $S$ for all $i\in I_1$. Since there are $2g$ numbers in $I_1\cup I_2$, we conclude that $S$ has $2g/2=g$ gaps. Thus, $S$ is reflective.
\end{proof}


\begin{figure}
\hspace{-3in}
\scalebox{.55}
{\small
\centering
    
    \setlength{\unitlength}{.8\linewidth}%
	\mbox{
    \begin{picture}(1, .4)
      \divide\unitlength by 8\relax
      \textcolor{gray}{%
        \put(0, 0){\makerule{1}{1}}%
        \put(4.5, 0){\makerule{1}{1}}%
        \put(9, 0){\makerule{1}{1}}%
        \put(13.5, 0){\makerule{1}{1}}%
        \put(1, 1){\makerule{3.5}{1}}%
        \put(5.5, 1){\makerule{3.5}{1}}%
        \put(10, 1){\makerule{3.5}{1}}%
        \put(14.5, 1){\makerule{2}{1}}%
        \put(0, 2){\makerule{16.5}{1}}%
      }%
      \put(0,0){\makebox(1,1){$0$}}
      \put(1,0){\makebox(1,1){$1$}}
      \put(2,0){\makebox(1.5,1){$\dots$}}
      \put(3.5,0){\makebox(1,1){$a-1$}}
      \put(4.5,0){\makebox(1,1){\circled{$a$}}}
      \put(6.5,0){\makebox(1.5,1){$\dots$}}
      \put(9,0){\makebox(1,1){$2a$}}
      \put(10,0){\makebox(1,1){$\dots$}}
      \put(11,0){\makebox(1.5,1){$\dots$}}
      \put(12.5,0){\makebox(1,1){$\dots$}}
      \put(13.5,0){\makebox(1,1){$qa$}}
      \put(14.5,0){\makebox(1,1){$\dots$}}
      \put(15.5,0){\makebox(1,1){$g-1$}}
      \put(0,1){\makebox(1,1){$g$}}
      \put(1,1){\makebox(1,1){\circled{$g+1$}}}
      \put(2,1){\makebox(1.5,1){\circled{$\dots$}}}
      \put(3.5,1){\makebox(1,1){\circled{$g+a-1$}}}
      \put(4.5,1){\makebox(1,1){$g+a$}}
      \put(6.5,1){\makebox(1.5,1){$\dots$}}
      \put(9,1){\makebox(1,1){$g+2a$}}
      \put(10,1){\makebox(1,1){$\dots$}}
      \put(11,1){\makebox(1.5,1){$\dots$}}
      \put(12.5,1){\makebox(1,1){$\dots$}}
      \put(13.5,1){\makebox(1,1){$g+qa$}}
      \put(14.5,1){\makebox(1,1){$\dots$}}
      \put(15.5,1){\makebox(1,1){$2g-1$}}\put(0,1){\makebox(1,1){$g$}}
      \put(0,2){\makebox(1,1){$2g$}}
      \put(1,2){\makebox(1,1){$\dots$}}
      \put(2,2){\makebox(1.5,1){\circled{$2g+a-r$}}}
      \put(3.5,2){\makebox(1,1){$\dots$}}
      \put(4.5,2){\makebox(1,1){$2g+a$}}
      \put(6,2){\makebox(1,1){$\dots$}}
      \put(9,2){\makebox(1,1){$2g+2a$}}
      \put(10,2){\makebox(1,1){$\dots$}}
      \put(11,2){\makebox(1.5,1){$\dots$}}
      \put(12.5,2){\makebox(1,1){$\dots$}}
      \put(13.5,2){\makebox(1,1){$2g+qa$}}
      \put(14.5,2){\makebox(1,1){$\dots$}}
      \put(15.5,2){\makebox(1,1){$3g-1$}}
      \multiput(0, 0)(0, 1){4}{\line(1, 0){16.5}}
      \multiput(0, 0)(1, 0){3}{\line(0, 1){3}}
      \multiput(3.5, 0)(1, 0){4}{\line(0, 1){3}}
      \multiput(8, 0)(1, 0){4}{\line(0, 1){3}}
      \multiput(12.5, 0)(1, 0){5}{\line(0, 1){3}}
      
    \end{picture}%
    }}
    \caption{For $g\in\N$, $a\in\ints{2,g+1}$ with $a\nmid g$, and integers $q,r$ for which \mbox{$g=qa+r$} with $r\in\ints{1,a-1}$, visualization of the reflective numerical semigroup $S=\left\langle \{a,2g+a-r\}\cup\{g+1, g+2, \dots, g+(a-1)\}\right\rangle$. Elements of $S$ are shaded. Generating elements of $S$ are circled. (We visualize $S$ within $\ints{0,3g-1}$ so that we can see $2g+a-r$.)
    }
    \label{fig:reflective-picture}
\end{figure}

See Figure~\ref{fig:reflective-picture} for a visualization of a general reflective numerical semigroup of genus $g$ and multiplicity $a$ as described in Proposition~\ref{prop:ED iff cases}.

\begin{example}\label{ex:count-reflective-genus-6}
    Suppose we want to describe all of the reflective numerical semigroups of genus $g=6$. By Proposition~\ref{prop:ED iff cases}, if $S$ is such a numerical semigroup with $\genus(S)=6$, then $\mult(S)\in\ints{2,7}$ with $\mult(S)\nmid 6$. Thus, $\mult(S)=a$ for $a\in\{4,5,7\}$. In what follows, recall that $r\equiv g\pmod*{a}$.
    
    If $a=4$, then $r\equiv 6\equiv2\pmod*{4}$. We have the numerical semigroup $S=\langle \{4,14\}\cup\ints{7,9}\rangle$. Then $\HH(S)=\{1,2,3,5,6,10\}$, which is equidistributed modulo 6. In this case, $S=\langle A\rangle$ is reflective for the minimal generating set $A=\{4,7,9\}$. Thus, $\embdim(S)=3$.
    
    If $a=5$, then $r\equiv 6\equiv1\pmod*{5}$. We have the numerical semigroup $S=\langle \{5,16\}\cup\ints{7,10}\rangle$. Then $\HH(S)=\{1,2,3,4,6,11\}$, which is equidistributed modulo 6. In this case, $S=\langle A\rangle$ is reflective for the minimal generating set $A=\{5,7,8,9\}$. Thus, $\embdim(S)=4$.
    
    If $a=7$, then $r\equiv 6\pmod*{7}$. We have the numerical semigroup $S=\langle \{7,13\}\cup\ints{7,12}\rangle$. Then $\HH(S)=\{1,2,3,4,5,6\}$, which is equidistributed modulo 6. In this case that $S=\langle A\rangle$ is reflective for the minimal generating set $A=\{7,8,9,10,11,12,13\}$. Thus, $\embdim(S)=7$.
    
    These are the three reflective numerical semigroups of genus 6. In Section~\ref{sec:counting-ME-semigroups}, we will look more closely at the problem of counting the number of reflective numerical semigroups of a given genus or of a given Frobenius number.
\end{example}

\subsection{Properties of reflective numerical semigroups}\label{subsec:ME-props}
For the remainder of this paper, we will assume $g\ge1$ and $a\in\ints{2,g+1}$ with $a\nmid g$. We therefore have $q,r\in\N_0$ such that $g=qa+r$ with $0<r<a$. In particular, $r$ is the unique integer in $\ints{1,a-1}$ for which $g\equiv r\pmod*{a}$.

In the proof of Proposition \ref{prop:ED iff cases}, we obtained a formula for the Frobenius number of a reflective numerical semigroup. We record it here.
\begin{cor}\label{cor:ME-frobenius}
    If $S$ is a reflective numerical semigroup with $\genus(S)=g\ge1$ and $\mult(S)=a$, then $\frob(S)=2g-r$.
\end{cor}

The set $A$ given in Proposition \ref{prop:ED iff cases} is, in general, not minimal because $g+(a-r)\in A$ and $a\mid \left(g+(a-r)\right)$. The form for the minimal generating set of $S$ depends on the congruence class of $g$ modulo $a$.

\begin{prop}\label{prop:min gen sets}
    Let $S$ be a reflective numerical semigroup of genus $\genus(S)=g\ge1$ and multiplicity $\mult(S)=a$. Then $S=\langle A\rangle$ for the minimal generating set $A$ given by
    \[A=
        \begin{dcases*} 
        \{a\}\cup\ints{g+1,g+(a-1)}\setminus\{g+(a-r)\} & if $g\not\equiv-1\pmod*{a}$, \\
        \{a,2g+1\}\cup\ints{g+2,g+(a-1)} & if $g\equiv-1\pmod*{a}$.
        \end{dcases*}
    \]
\end{prop}

\begin{proof}
    By Proposition \ref{prop:ED iff cases}, if $S$ is reflective, then $S=\langle B\rangle$ for the generating set 
    \[B=\{a,2g+a-r\}\cup\ints{g+1,g+(a-1)}.\] 
    To obtain the minimal generating set $A$ for $S$, we need to eliminate elements of $B$ which are non-negative integer combinations of smaller elements of $B$. Once we have eliminated all such elements, we will have a minimal generating set $A$. In what follows, recall that $g\equiv r\pmod*{a}$.
    
    To start, we first observe that the minimal element of $B$ cannot be written as a non-negative linear combination of other elements of $B$. Thus, $a=\mult(S)$ must be in $A$.
    
    Next, we consider $2g+a-r\in B$. We have two cases for $r$: either $r\in\ints{1,a-2}$; or $r=a-1$.
    
    If $r\in\ints{1,a-2}$, then $2g+a-r=(g+1)+(g+a-r-1)$ with $g+1,g+a-r-1\in B$. Since $2g+a-r$ is a non-negative linear combination of other elements of $B$, we have $2g+a-r\not\in A$.
    
    If $r=a-1$, then $g+(q+1)a=2g+1$. By Lemma \ref{lem:toward-minimality}, if $2g+1$ is a non-negative linear combination of other elements of $B$, then $2g+1=la+(g+i)$ for some $l\in\mathbb{N}_0$ and  $i\in\ints{1,a-1}$. However, $2g+1\equiv-1\pmod*{a}$ and $la+(g+i)\equiv i-1\pmod*{a}$ for some $i\in\ints{1,a-1}$. Thus $i\equiv0\pmod*{a}$ and $i\in\ints{1,a-1}$, a contradiction. Therefore, when $r=a-1$, $2g+a-r=2g+1$ is not a non-negative linear combination of other elements of $B$. In this case, we have $2g+a-r=2g+1\in A$.
    
    Finally, consider $g+j\in B$ for some $j\in\ints{1,a-1}$. If $g+j$ is a non-negative linear combination of other elements of $B$, then, again by Lemma \ref{lem:toward-minimality}, we may assume that that combination contains at most one element of the form $g+i$ for $i\in\ints{1,a-1}$. We therefore have one of the following: $g+j=la$ for $l\in\N_0$; or $g+j=la+(g+i)$ for $l\in\N_0$ and $i\in\ints{1,a-1}$.
    
    The first case occurs precisely when $a\mid g+j$, which occurs when $j=a-r$. Thus, $g+(a-r)\not\in A$. 
    
    For the second case, if we consider the equation $g+j=la+(g+i)$ modulo $a$ we see that $i\equiv j\pmod*{a}$. This means $i=j$ and thus $l=0$. We can therefore only write $g+j$ in terms of itself. From our analysis of the first case, we see that the second case occurs precisely when $j\ne a-r$. When this occurs, we cannot write $g+j$ as a non-negative linear combination of other elements of $B$. We conclude that $g+j\in A$ for all $j\in\ints{1,a-1}\setminus\{a-r\}$.
    
    Putting everything together, starting with the generating set $B$ from Proposition \ref{prop:ED iff cases} for a reflective numerical semigroup $S$, we eliminate $g+(a-r)$. In the case where $g\not\equiv-1\pmod*{a}$, we also eliminate $2g+a-r$. The result is a minimal generating set $A$ which generates the same reflective numerical semigroup $S$.
\end{proof}

From Proposition \ref{prop:min gen sets}, we have the minimal generating set for a reflective numerical semigroup. We record the embedding dimension of such a numerical semigroup in the corollary below.

\begin{cor}\label{cor:ME-dimension}
    Let $S$ be a reflective numerical semigroup of genus $\genus(S)=g\ge1$ and multiplicity $\mult(S)=a$. Then the embedding dimension of $S$ is given by 
    \[
        \embdim(S) = 
        \begin{dcases*} 
            a-1 & if $g\not\equiv-1\pmod*{a}$, \\ a & if $g\equiv-1\pmod*{a}$.
        \end{dcases*}
    \]
\end{cor}

We are almost ready to explicitly compute an Ap\'ery set of a reflective numerical semigroup. To do so, we first need a lemma.

\begin{lem}\label{lem:subset-of-apery-set}
    Suppose $S=\langle A\rangle$ with $A=\{a_1,\dots,a_k\}$ minimal. If $a_i\not\equiv a_j\pmod*{a_1}$ for all $i,j\in\ints{2,k}$ with $i\ne j$, then \[\{0,a_2,\dots,a_k\}\subseteq\Ap(S;a_1).\]
    Furthermore, if $\#A=a_1$, then $\{0,a_2,\dots,a_k\}=\Ap(S;a_1)$.
\end{lem}
\begin{proof}
    Clearly, $0\in\Ap(S;a_1)$.

    Mimicking the argument given in \cite[Page 20, second paragraph]{RosalesGarciaSanchez09}, since $a_i$ is in the minimal generating set for $S$, $a_i-x\not\in S$ for all $x\in S\setminus\{0,a_i\}$. In particular, $a_i-a_1\not\in S$ for all $i\ge2$. Since $a_i\in S$, this means $a_i\in\Ap(S;a_1)$.
    
    Therefore, $\{0,a_2,\dots,a_k\}\subseteq\Ap(S;a_1)$. And since $\#\Ap(S;a_1)=a_1$, if $k=a_1$ then we have equality.
\end{proof}

\begin{cor}\label{cor:ME_Apery_set}
    Let $S$ be a reflective numerical semigroup with $\genus(S)=g\ge1$ and $\mult(S)=a$. Then the Ap\'ery set of $S$ relative to $a$ is
    \[\Ap(S;a) = \{0,2g+a-r\}\cup\ints{g+1,g+(a-1)}\setminus\{g+(a-r)\}.\]
    When $g\equiv-1\pmod*{a}$, we can write this more compactly as \[\Ap(S;a)=\{0,2g+1\}\cup\ints{g+2,g+(a-1)}.\]
\end{cor}
\begin{proof}
    Suppose $g\equiv-1\pmod*{a}$. By Proposition~\ref{prop:min gen sets}, $S$ is minimally generated by 
    \[
        A=\{a,2g+1\}\cup\ints{g+2,g+(a-1)}.
    \]
    Since $a\mid(g+1)$, $g+i\equiv (i-1)\pmod*{a}$ for each $i\in\ints{2,a-1}$, and $2g+1\equiv-1\pmod*{a}$. Thus $A$ contains one element in each congruence class modulo $a$. By Lemma \ref{lem:subset-of-apery-set},  
    \[
    \Ap(S;a)=\{0,2g+1\}\cup\ints{g+2,g+(a-1)}.
    \]
    To write this in terms of $r$, note that $g\equiv-1\pmod*{a}$ implies $r=a-1$.  We can therefore write the Ap\'ery set as
    \[
        \Ap(S;a)=\{0,2g+a-r\}\cup\ints{g+1,g+(a-1)}\setminus\{g+(a-r)\}.
    \]

    Now suppose $g\not\equiv-1\pmod*{a}$. By Proposition~\ref{prop:min gen sets}, $S$ is minimally generated by 
    \[
        A=\{a\}\cup\ints{g+1,g+(a-1)}\setminus\{g+(a-r)\}.
    \] 
    The $a-2$ elements of $A\setminus\{a\}$ are all distinct modulo $a$. By Lemma \ref{lem:subset-of-apery-set}, $(\{0\}\cup A\setminus\{a\})\subseteq\Ap(S;a)$. We therefore have elements in $\Ap(S;a)$ from every congruence class except for $g$ modulo $a$. To identify the final element, since $\frob(S)\equiv g\pmod*{a}$ and $\frob(S)=\max\left(\HH(S)\right)$, the smallest element of $S$ that is congruent to $g$ modulo $a$ is $\frob(S)+a=2g-r+a$. Therefore, $2g-r+a\in\Ap(S;a)$ as well.
\end{proof}

If we have the Ap\'ery set of a numerical semigroup $S$ relative to its multiplicity, then we have a simple criterion to determine if $S$ is reflective or not. This follows immediately from  Corollary~\ref{cor:ME_Apery_set} and Proposition~\ref{prop:same-apery-same-semigroup}.

\begin{cor}
    Suppose $S$ is a numerical semigroup with $\genus(S)\ge1$, $\mult(S)\in\ints{2,\genus(S)+1}$, and $\mult(S)\nmid\genus(S)$. For $r\in\ints{1,\mult(S)}$ with $r\equiv \genus(S)\pmod*{\mult(S)}$, we have that $S$ is reflective if and only if
    \[\Ap\left(S;\mult(S)\right)=\{0,2\genus(S)+\mult(S)-r\}\cup\ints{\genus(S)+1,\genus(S)+\mult(S)-1}\setminus\{\genus(S)+\mult(S)-r\}.\]
\end{cor}

\section{Characterizations of reflective numerical semigroups}\label{sec:ME-props}
In this section, we'll describe a few commonly-studied families of numerical semigroups and determine which members of those families are reflective.

\subsection{Symmetric, pseudo-symmetric, and irreducible numerical semigroups}
We begin with symmetric, pseudo-symmetric, and irreducible numerical semigroups. For a thorough treatment, see \cite[Chapter 3]{RosalesGarciaSanchez09}.

As mentioned in the introduction, a numerical semigroup $S$ is \emph{symmetric} if exactly one of $z$ and $\frob(S)-z$ is in $S$ for all $z\in\Z$. A numerical semigroup $S$ is \emph{pseudo-symmetric} if $\frob(S)$ is even and exactly one of $z$ and $\frob(S)-z$ is in $S$ for all $z\in\Z\setminus\{\frob(S)/2\}$. A numerical semigroup $S$ is \emph{irreducible} if $S$ is not the intersection of two numerical semigroups which properly contain $S$.

\begin{prop}[{\cite[Corollary 4.5]{RosalesGarciaSanchez09}}]
\label{prop:symmetric-criterion}
    A numerical semigroup $S$ is symmetric if and only if $\frob(S)=2\genus(S)-1$.
\end{prop}

\begin{prop}[{\cite[Lemma 2.2]{Garcia-MarcoRamirez-AlfonsinRodseth2017}}]
\label{prop:pseudo-symmetric-criterion}
    A numerical semigroup $S$ is pseudo-symmetric if and only if $\frob(S)=2\genus(S)-2$.
\end{prop}

\begin{prop}[{\cite[Proposition 2]{RosalesBranco2003}}]\label{prop:irred-criterion}
    A numerical semigroup $S$ is irreducible if and only if $S$ is either symmetric or pseudo-symmetric.
\end{prop}

The set of \emph{pseudo-Frobenius numbers} of $S$ is
\[
    \PF(S) 
    = \{h\in \HH(S) : h+s\in S\text{ for all } 0\ne s\in S\}.
\]
This is the set of gaps for which adding any nonzero semigroup element results in a semigroup element. Since $\frob(S)$ is the largest gap, we clearly have $\frob(S)\in\PF(S)$. The \emph{type} of $S$ is $\type(S)=\#\PF(S)\ge1$, the number of pseudo-Frobenius numbers of $S$. Since every element of $S$ can be written as a finite sum of elements from a generating set of $S$, we obtain the following characterization of $\PF(S)$.
\begin{lem}\label{lem:alternative-pseudo-frob}
    Suppose $S=\langle A\rangle$ for $A=\{a_1,a_2,\dots,a_k\}\subset\N$. Then 
    \[
        \PF(S)=\left\{h\in\HH(S) : h+a_i\in S\text{ for all } i\in\ints{1,k}\right\}.
    \]
\end{lem}
\begin{proof}
    Let $Y=\{h\in\HH(S) : h+s\in S\text{ for all } 0\ne s\in S\}$ and $Z = \{h\in\HH(S) : h+a_i\in S\text{ for all } i\in\ints{1,k}\}.$ We wish to show $Y=Z$. We will do so by showing $Y\subseteq Z$ and $Z\subseteq Y$.
    
    Since $0\ne a_i\in S$ for all $i$, we immediately see $Y\subseteq Z$.
    
    We now wish to show $Z\subseteq Y$. Suppose $h\in Z$. Then $h\in\HH(S)$ with $h+a_i\in S$ for all $i\in\ints{1,k}$. For any nonzero $s\in S$, since $S=\langle a_1,a_2,\dots,a_k\rangle$, there are non-negative integers $b_1,b_2,\dots,b_k$ such that 
    \[s=c_1a_1+c_2a_2+\dots+c_ka_k.\]
    Since $s\ne0$, there is some $l\in\ints{1,k}$ for which $b_l>0$. We obtain non-negative integers $c_1,c_2,\dots,c_k$ for which \[s=a_l + \left(c_1a_1+c_2a_2+\dots+c_ka_k\right)\]
    by letting $c_l=b_l-1$ and, for $j\ne l$, $c_j=b_j$. In particular, $c_1,c_2,\dots,c_k\ge0$. We then have
    \[h+s=(h+a_l) + \left(c_1a_1+c_2a_2+\dots+a_kc_k\right).\]
    By assumption, $h+a_l\in S$. Hence we have $h+s$ written as a sum of two elements of $S$. We therefore have $h+s\in S$. Since $s$ is an arbitrary nonzero element of $S$, we must have $h\in Y$. Thus, $Z\subseteq Y$.
\end{proof}

\begin{prop}[{\cite[Corollary 4.11 \& Corollary 4.16]{RosalesGarciaSanchez09}}]\label{prop:symmetric-and-pseudo-symmetric-criterioa}
    A numerical semigroup $S$ is symmetric if and only if $\type(S)=1$. A numerical semigroup $S$ is pseudo-symmetric if and only if $\PF(S)=\{\frob(S)/2,\frob(S)\}$.
\end{prop}

We can now compute the pseudo-Frobenius numbers of a reflective numerical semigroup.
\begin{prop}\label{prop:PF(S)}
    Let $S$ be a reflective numerical semigroup with $\genus(S)=g\ge1$ and $\mult(S)=a$. Let $q=\Floor{g/a}$ and $r=g-qa$. Then 
    \[\PF(S) = \ints{g-(r-1),g-1}\cup\{2g-r\}.\]
    In particular, $\type(S)=r$. 
\end{prop}
\begin{proof}
    Since $\PF(S)\subset\HH(S)$, we first partition $\HH(S)$ into three disjoint sets 
    \[\HH(S)=H_1\sqcup H_2\sqcup H_3,\]
    where $H_1 = \{h\in \HH(S) : h\le qa\}$, $H_2 = \{h\in\HH(S) : qa<h<g\}$, $H_3 = \{h\in\HH(S) : g\le h\}$.
    
    Let $h\in H_1$. Since $H_1\subset\HH(S)$ and $a\in S$, we have $a\nmid h$. Thus, $qa-h=la+i$ for some $l\in\N_0$ and $i\in\ints{1,a-1}$. Since $la\in S$ and $g+i\in S$, their sum $la+(g+i)$ is a positive element of $S$. Then $h+(la+(g+i))=2g-r=\frob(S)\not\in S$, and thus $h\not\in\PF(S)$.
    
    Let $h\in H_2$. Applying Lemma \ref{lem:alternative-pseudo-frob} and Proposition \ref{prop:ED iff cases}, we will show $h+\alpha\in S$ for all $\alpha\in A=\{a,2g+a-r\}\cup\ints{g+1,g+(a-1)}$. First, let $\alpha=a$. Since $g=qa+r$, we have $h=qa+j$ for some $j\in\ints{1,r-1}$. Then $g<h+a<g+a$, so $h+a=g+i$ for some $i\in\ints{1,a-1}$. Thus $h+\alpha=h+a=g+i\in S$. Next, let $\alpha=2g+a-r$. We have $h+\alpha=h+(2g+a-r)>\frob(S)$, so $h+\alpha\in S$. Finally, let $\alpha=g+i$ for some $i\in\ints{1,a-1}$. Then $h+\alpha=h+(g+i)=2g-r+(i+j)>\frob(S)$, so $h+\alpha\in S$. Since $h+\alpha\in S$ for all $h\in H_2$ and all $\alpha\in A$, $H_2\subset\PF(S)$.
    
    Let $h\in H_3$. Then $h=g+ja$ for some $j\in\ints{0,q}$. Clearly $h+a\not\in S$ for $j\in\ints{0,q-1}$. When $j=q$, $h=g+qa=2g-r=\frob(S)\in\PF(S)$. Thus, $H_3\cap\PF(S)=\{2g-r\}$.
    
    To conclude, $\PF(S)=H_2\cup\{2g-r\}=\ints{g-(r-1),g-1}\cup\{2g-r\}$, a set of $r$ integers.
\end{proof}

Now that we have $\type(S)$, we can determine conditions on $\mult(S)$ and $\genus(S)$ for which a reflective numerical semigroup $S$ is symmetric or pseudo-symmetric. 

\begin{cor}\label{cor:ME and symm or pseudo-symm}
    Let $S$ be a reflective numerical semigroup with $\genus(S)=g$ and $\mult(S)=a$. Then
    \begin{itemize}
        \item $S$ is symmetric if and only if $g\equiv1\pmod*{a}$;
        \item $S$ is pseudo-symmetric if and only if $g\equiv2\pmod*{a}$; and
        \item $S$ is irreducible if and only if $g\equiv1,2\pmod*{a}$.
    \end{itemize}
\end{cor}
\begin{proof}
    We use Proposition \ref{prop:symmetric-and-pseudo-symmetric-criterioa} and Proposition \ref{prop:PF(S)}. If $g\equiv1\pmod*{a}$, then $\type(S)=1$, which implies that $S$ is symmetric. If $g\equiv2\pmod*{a}$, then $\PF(S)=\{g-1,2g-2\}$. By Corollary~\ref{cor:ME-frobenius}, $\frob(S)=2g-2$. Hence $\PF(S)=\{\frob(S)/2,\frob(S)\}$, which implies that $S$ is pseudo-symmetric. If $g\not\equiv1,2\pmod*{a}$, then $\type(S)>2$, so $S$ is neither symmetric nor pseudo-symmetric. Finally, by Proposition \ref{prop:irred-criterion}, $S$ is irreducible if and only if $S$ is symmetric or pseudo-symmetric.
    \end{proof}

\begin{remark}
    From Corollary~\ref{cor:ME and symm or pseudo-symm} we see that a reflective numerical semigroup $S$ with $\mult(S)=2$ is necessarily symmetric. (As we will see in the next subsection, for any numerical semigroup $S$, if $\mult(S)=2$, then $\embdim(S)=2$, and this alone implies that $S$ is symmetric.) If $\mult(S)>2$ (and hence $g\equiv1\pmod*{a}$ implies $g\not\equiv-1\pmod*{a}$), we use Proposition \ref{prop:min gen sets} to see that a reflective numerical semigroup $S$ is symmetric if and only if $g=ka+1$ for some $k\in\N$. In this case, $S$ is generated by a set of the form
    \[A = \{a\}\cup\ints{g+1,g+(a-2)} = \{a\}\cup\ints{ka+2,ka+a-1}.\]
\end{remark}

\subsection{Numerical semigroups generated by generalized arithmetic sequences}
We now focus on numerical semigroups generated by generalized arithmetic sequences (defined below). Once we have determined the reflective members of this family, we can specialize to find all reflective numerical semigroups that are generated by an arithmetic sequence. This includes a description of all reflective numerical semigroups with embedding dimension 2 (i.e., $S=\langle a,b\rangle$).

For $a,h,d,k\in\N$, a \emph{generalized arithmetic sequence} is a sequence of length $k$ of the form
\[
    (a, ha+d, ha+2d, ha+3d, \dots, ha+(k-1)d).
\]
Numerous properties of numerical semigroups that are generated by generalized arithmetic sequences have been determined. See, e.g., \cite{Lewin1975}, \cite{Selmer1977}, \cite{Ritter1998}, and \cite{Matthews2004}. For reflective numerical semigroups generated by generalized arithmetic sequences in the case where $k=a$, see \cite[Corollary 5.6]{Shor22}.

We will now describe the reflective numerical semigroups that are generated by generalized arithmetic sequences. Once we have this result, we will obtain results for numerical semigroups generated by arithmetic sequences and for numerical semigroups with embedding dimension 2.

\begin{prop}\label{prop:ME-gen-by-generalized-arith-seq}
    Suppose $S=\langle a, ha+d, ha+2d, \dots, ha+(k-1)d\rangle$ is a reflective numerical semigroup with $\genus(S)=g$. Then $S$ belongs to one of the five following families.
    \begin{enumerate}
        \item $S=\langle 1\rangle=\N_0$. In this case, $g=0$, $a=1$, $k=1$, and $h$ and $d$ are arbitrary.
        \item $S=\langle 2, 2g+1\rangle$ for $g\equiv1\pmod*{2}$. In this case, $a=2$, $k=2$, and $2h+d=2g+1$.
        \item $S=\langle 3, g+1\rangle$ for $g\equiv1\pmod*{3}$. In this case, $a=3$, $k=2$, and $3h+d=g+1$.
        \item $S=\langle 3, g+2, 2g+1\rangle$ for $g\equiv2\pmod*{3}$. In this case, $a=3$, $k=3$, and \mbox{$d=g-1\equiv1\pmod*{3}$.}
        \item $S=\langle \ints{g+1,2g+1}\rangle$ for any $g\ge0$. In this case, $a=g+1$, $k=g+1$, and $d=1$.
    \end{enumerate}
\end{prop}
\begin{proof}
    Suppose $S$ is reflective and $\genus(S)=g$. By Prop.~\ref{prop:min gen sets}, we have the explicit form of a generating set of $S$. We consider cases based on $\mult(S)$.
    
    If $\mult(S)=1$, then $S=\langle 1\rangle=\N_0$, a numerical semigroup of genus $\genus(S)=0$. The sequence $(1)$ is a generalized arithmetic sequence where $a=1$, $k=1$, and $h$ and $d$ are arbitrary.
    
    If $\mult(S)=2$, then by Prop.~\ref{prop:min gen sets}, $S=\langle 2, 2g+1\rangle$ for odd $g$. The sequence $(2,2g+1)$ is a generalized arithmetic sequence where $a=2$, $k=2$, and $2h+d=2g+1\equiv1\pmod*{2}$.
    
    If $\mult(S)=3$, then by Prop.~\ref{prop:min gen sets}, we have either $S=\langle 3, g+1\rangle$ for $g\equiv1\pmod*{3}$, or $S=\langle 3, g+2, 2g+1\rangle$ for $g\equiv2\pmod*{3}$. For the first case, the sequence $(3,g+1)$ is a generalized arithmetic sequence where $a=3$, $k=2$, and $3h+d=g+1$. For the second case, the sequence $(3,g+2,2g+1)$ is a generalized arithmetic sequence where $a=3$, $k=3$, $d=g-1$, and $h=1$.
    
    If $\mult(S)\ge4$, we consider three cases: $\mult(S)=4$ and $g\equiv2\pmod*{4}$; $\mult(S)\ge4$ and $g\not\equiv-1\pmod*{4}$, except for the combination of $\mult(S)=4$ and $g\equiv2\pmod*{4}$; and $\mult(S)\ge4$ and $g\equiv-1\pmod*{4}$.
    
    For the first case, by Prop.~\ref{prop:min gen sets}, we have $S=\langle 4, g+1, g+3\rangle$. Then $d=(g+3)-(g+1)=2$ and $g+1=4h+d=4h+2$. This implies $g\equiv1\pmod*{4}$. However, we assumed that $g\equiv2\pmod*{4}$ so this is not possible.
    
    For the second case, by Prop.~\ref{prop:min gen sets}, for $a=\mult(S)\ge4$, the set $\ints{g+1,g+(a-1)}\setminus\{g+(a-r)\}$ contains at least two consecutive integers. (This is why we have specifically excluded the combination of $a-1=3$ and $r=2$ here.) Thus, $d=1$. Since $r\ne a-1$, $g+1$ is one of the generating elements. Thus, $ha+d=g+1$. However, since $d=1$, this implies $a\mid g$, a contradiction. This is not possible.
    
    For the third case, by Prop.~\ref{prop:min gen sets}, for $a=\mult(S)\ge4$ we have $S=\langle A\rangle$ for $A=\{a,2g+1\}\cup\ints{g+2,g+(a-1)}$. $A$ contains consecutive integers $g+2$ and $g+3$. Therefore, $d=1$.  The two largest elements of this generating set are $g+(a-1)$ and $2g+1$. Since their difference is $d$, which is 1, we must have $g=a-1$. It follows that $A=\ints{g+1,2g+1}$.
    
    In the previous paragraph, we assumed that $a\ge4$ and hence $g=a-1\ge3$. However, we note that $S=\langle \ints{g+1,2g+1}\rangle$ is reflective for $g=0,1,2$ as well. These numerical semigroups are special cases of the numerical semigroups described in cases (1), (2), and (4) in the statement of the proposition.
\end{proof}

We immediately obtain two corollaries for reflective numerical semigroups -- those that are generated by arithmetic sequences, and those with embedding dimension 2.

An arithmetic sequence is a generalized arithmetic sequence with $h=1$. With that in mind, we have the following.

\begin{cor}\label{cor:ME-gen-by-arithmetic-seq}
    Suppose $S=\langle a, a+d, a+2d, \dots, a+(k-1)d\rangle$ for $a,d,k\in\N$ with $k\le a$, a numerical semigroup generated by an arithmetic sequence of length $k$. Then $S$ is reflective if and only if one of the following cases occurs.
    \begin{enumerate}
        \item $S=\langle 1\rangle$. In this case, $g=0$, $a=1$, $k=1$, and $d$ is arbitrary.
        \item $S=\langle 2, 2g+1\rangle$ for $g\equiv1\pmod*{2}$. In this case, $a=2$, $k=2$, and $d=2g-1$.
        \item $S=\langle 3, g+1\rangle$ for $g\equiv1\pmod*{3}$. In this case, $a=3$, $k=2$, and $d=g-2$.
        \item $S=\langle 3, g+2, 2g+1\rangle$ for $g\equiv2\pmod*{3}$. In this case, $a=3$, $k=3$, and $d=g-1$.
        \item $S=\langle \ints{g+1,2g+1}\rangle$ for any $g\ge0$. In this case, $a=g+1$, $k=g+1$, and $d=1$.
    \end{enumerate}
\end{cor}

Next, a numerical semigroup with embedding dimension 2 is a numerical semigroup of the form $S=\langle a,b\rangle$ with $1<a<b$ and $\gcd(a,b)=1$. We can view this as a numerical semigroup generated by a generalized arithmetic sequence of length 2. We use cases (2) and (3) from Prop.~\ref{prop:ME-gen-by-generalized-arith-seq}.

\begin{cor}\label{cor:ED dim2}
    Suppose $S$ is a numerical semigroup with $\embdim(S)=2$. Then $S$ is reflective if and only if for some $n\in\N_0$, $S=\langle A\rangle$ for $A=\{2,4n+3\}$ (in which case $g=2n+1$) or $A=\{3,3n+2\}$ (in which case $g=3n+1$).
\end{cor}

\subsection{Free numerical semigroups}
We now briefly describe free numerical semigroups, which were named in \cite{BertinCarbonne1977}.

\begin{defn}
    Given a sequence $A=(a_1,\dots,a_k)$ of non-negative integers, for $i\in\ints{1,k}$, let $A_i = (a_1,\dots,a_i)$, $d_i=\gcd(A_i)$, and $S_i=\langle A_i\rangle$. Let $\cc(A)=(c_2,\dots,c_k)$ where $c_j=d_{j-1}/d_j$ for $j\in\ints{2,k}$. If we have $c_ja_j\in S_{j-1}$ for all $j\in\ints{2,k}$, then we say $A=(a_1,\dots,a_k)$ is a \emph{telescopic sequence} (or a \emph{smooth sequence}).
\end{defn}
\begin{defn}
    If $A$ is a set which can be ordered to form a telescopic sequence, then we say $A$ is a \emph{telescopic set}. If $S$ is a numerical semigroup which is generated by a telescopic set $A$ (which necessarily has $\gcd(A)=1$), then we say $S$ is a \emph{free numerical semigroup}.
\end{defn}

For a free numerical semigroup $S$, we can construct an Ap\'ery set and compute the genus and Frobenius number of $S$.

\begin{prop}\label{prop:free-num-sgps-results}
    Let $S$ be a free numerical semigroup, so $S=\langle A\rangle$ for some telescopic sequence $A=(a_1,\dots,a_k)$ with $\cc(A)=(c_2,\dots,c_k)$. Then 
    \[
        \Ap(S;a_1) = \left\{\sum\limits_{j=2}^{k} n_j a_j : 0\le n_j < c_j \right\},
    \]
    \[
        \genus(S)=\dfrac{1}{2}\left(1-a_1+\sum\limits_{j=2}^k c_j(g_j-1)\right),
    \]
    and $\frob(S)=2\genus(S)-1$. In particular, $S$ is a symmetric numerical semigroup.
\end{prop}

Note that each element of the Ap\'ery set described in Prop.~\ref{prop:free-num-sgps-results} is the $(k-1)$-fold sum of terms from arithmetic sequences that each start at 0. If $S$ is a numerical semigroup of embedding dimension 2, then $S=\langle a,b\rangle$ with $\gcd(a,b)=1$, and $\Ap(S;a)=\{0,b,2b,\dots,(a-1)b\}$, an arithmetic sequence. As a result, we can see that a numerical semigroup of embedding dimension 2 is necessarily free. The similar structure of the Ap\'ery sets means that free numerical semigroups are something of a natural generalization of numerical semigroups of embedding dimension 2. (In particular, one can compute invariants of free numerical semigroups in a similar manner as we do for numerical semigroups of embedding dimension 2. See \cite[Theorem 2.3, Corollary 2.6]{GassertShor17}.)

A numerical semigroup is free if it is generated by a telescopic sequence. However, the corresponding telescopic set need not be minimal. Fortunately, the elements in the minimal generating set for the numerical semigroup can be ordered to form a telescopic sequence as was shown in \cite{Shor19a}.
\begin{prop}[{\cite[Theorem 52]{Shor19a}}]
\label{prop:can-assume-telescopic-set-minimal}
    Suppose $S$ is a free numerical semigroup and $S=\langle A\rangle$ for a minimal set $A$. Then $A$ is a telescopic set.
\end{prop}

Now we can characterize the reflective numerical semigroups which are free. For this characterization, we first need a few lemmas about free numerical semigroups.

\begin{lem}\label{lem:gcd-1-goes-last}
    Suppose $A=(a_1,\dots,a_k)$ is a minimal telescopic sequence. If $\gcd(a_i,a_j)=1$ for some $i<j$, then $j=k$.
\end{lem}
\begin{proof}
    For $A$ a minimal telescopic sequence with $\gcd(a_i,a_j)=1$ for some $i<j$, suppose $j\ne k$. Since $i<j$, $d_j=\gcd(a_1,a_2,\dots,a_j)=1$ and therefore $d_{j+1}=\gcd(a_1,\dots,a_{j+1})=1$ as well. Then $c_{j+1}=d_{j}/d_{j+1}=1/1=1$. Since $A$ is telescopic, $c_{j+1}a_{j+1}\in\langle a_1,\dots,a_j\rangle$. Thus $a_{j+1}\in\langle a_1,\dots,a_j\rangle$. But this means $a_{j+1}$ can be written as a non-negative linear combination of other elements of $A$, which contradicts the minimality of $A$. Hence $j=k$, as desired.
\end{proof}

\begin{lem}\label{lem:three-mutually-coprime}
    Suppose $A=(a_1,\dots,a_k)$ is a minimal telescopic sequence. Then for all distinct $x,y,z\in\ints{1,k}$, either $\gcd(a_x,a_y)>1$ or $\gcd(a_x,a_z)>1$ or $\gcd(a_y,a_z)>1$.
\end{lem}
\begin{proof}
    Suppose $A$ is a minimal telescopic sequence and let $x,y,z$ be distinct elements of $\ints{1,k}$. We have $\gcd(a_x,a_y)\ge1$, $\gcd(a_x,a_z)\ge1$, and $\gcd(a_y,a_z)\ge1$. If $\gcd(a_x,a_y)>1$ or $\gcd(a_x,a_z)>1$, then we are done. If $\gcd(a_x,a_y)=1$ and $\gcd(a_x,a_z)=1$, we apply Lemma~\ref{lem:gcd-1-goes-last}. Since $\gcd(a_x,a_y)=1$, one of $x$ and $y$ is equal to $k$. Since $\gcd(a_x,a_z)=1$, one of $x$ and $z$ is equal to $k$. We have that $x,y,z$ are distinct. Therefore we must have $x=k$, $y\ne k$, and $z\ne k$. By the contrapositive to Lemma~\ref{lem:gcd-1-goes-last}, we have $\gcd(a_y,a_z)\ne1$, from which we conclude that $\gcd(a_y,a_z)>1$, as desired.
\end{proof}

\begin{lem}\label{lem:distinct-coprime-pairs}
    Suppose $A=(a_1,\dots,a_k)$ is a minimal telescopic sequence. Then, for all distinct $x,y,z,w$ in $\ints{1,k}$, either $\gcd(a_x,a_y)>1$ or $\gcd(a_z,a_w)>1$.
\end{lem}
\begin{proof}
    Suppose $A$ is a minimal telescopic sequence and let $x,y,z,w$ be distinct elements of $\ints{1,k}$. We have $\gcd(a_x,a_y)\ge1$ and $\gcd(a_z,a_w)\ge1$. If $\gcd(a_x,a_y)>1$, then we are done. If $\gcd(a_x,a_y)=1$, then by Lemma~\ref{lem:gcd-1-goes-last}, either $x=k$ or $y=k$. Since $x,y,z,w$ are distinct, we have $z\ne k$ and $w\ne k$. By the contrapositive of Lemma~\ref{lem:gcd-1-goes-last}, we therefore have $\gcd(a_z,a_w)\ne1$, from which we conclude $\gcd(a_z,a_w)>1$, as desired.
\end{proof}

The main takeaway of Lemma~\ref{lem:three-mutually-coprime} and Lemma~\ref{lem:distinct-coprime-pairs} is that if a minimal set contains three elements which are pairwise relatively prime, or contains two distinct pairs of relatively prime integers, then the set is not a telescopic set.

We can now characterize all free numerical semigroups which are reflective. Recall that a free numerical semigroup is one which is generated by a telescopic set, and we may assume that such a telescopic set is minimal. (See Prop.~\ref{prop:can-assume-telescopic-set-minimal}.)

\begin{prop}\label{prop:ME-and-free}
    Suppose $S$ is a free numerical semigroup that is reflective. Then $S=\langle A\rangle$ for a minimal telescopic set $A$ given by one of the following.
    \begin{enumerate}
        \item $A=\{1\}$ (and thus $\genus(S)=0$).
        \item $A=\{2,4n+3\}$ for $n\in\N_0$ (and thus $\genus(S)=2n+1$).
        \item $A=\{3,3n+2\}$ for $n\in\N_0$ (and thus $\genus(S)=3n+1$).
        \item $A=\{4,4n+2,4n+3\}$ for $n\in\N_0$ (and thus $\genus(S)=4n+1$).
    \end{enumerate}
\end{prop}
\begin{proof}
To begin, we assume that $S$ is reflective. Then $S=\langle A\rangle$ for a minimal set $A$ as is given in Prop.~\ref{prop:min gen sets}. We want to determine when $A$ is a telescopic set, which means $S$ is a free numerical semigroup. We will consider cases based on $\mult(S)$. Since $S$ is reflective, if $g\ge1$, then $\mult(S)\nmid g$.

If $\mult(S)=1$, then $S=\langle A\rangle$ for $A=\{1\}$, a telescopic set. Thus $S=\N_0$ is free and $\genus(S)=0$.

If $\mult(S)=2$, then $S=\langle 2, b\rangle$ for some odd $b\ge3$. Note that $(2,b)$ is a telescopic sequence for any odd $b$. Since $\embdim(S)=2$, we use Corollary~\ref{cor:ED dim2} to conclude that $S$ is reflective precisely when $b=4n+3$ for $n\in\N_0$. In this case, $\genus(S)=2n+1$.

If $\mult(S)=3$, then $g\equiv1,2\pmod*{3}$. We will consider these cases separately.

Suppose $\mult(S)=3$ and $g\equiv1\pmod*{3}$. Then $g=3n+1$ for some $n\in\N_0$, and $S=\langle 3,g+1\rangle=\langle 3,3n+2\rangle$. Note that $(3,3n+2)$ is a telescopic sequence, so $S=\langle 3,3n+2\rangle$ is a free numerical semigroup.

Suppose $\mult(S)=3$ and $g\equiv2\pmod*{3}$. Then $g=3n+2$ for some $n\in\N_0$, and $S=\langle 3,g+2,2g+1\rangle=\langle 3, 3n+4, 6n+5\rangle$. Since $3$, $3n+4$, and $6n+5$ are pairwise relatively prime, we apply Lemma~\ref{lem:three-mutually-coprime} to conclude that $S$ is not a free numerical semigroup.

If $\mult(S)=4$, then $g\equiv1,2,3\pmod*{4}$. We will consider these cases separately.

Suppose $\mult(S)=4$ and $g\equiv1\pmod*{4}$. Then $g=4n+1$ for some $n\in\N_0$, and $S=\langle 4, g+1, g+2\rangle=\langle 4, 4n+2, 4n+3\rangle$. Then $c_2=c_3=2$, and one confirms that $(4,4n+2,4n+3)$ is a telescopic sequence. Hence, $S=\langle 4,4n+2,4n+3\rangle$ is a free numerical semigroup.

Suppose $\mult(S)=4$ and $g\equiv2\pmod*{4}$. Then $g=4n+2$ for some $n\in\N_0$, and $S=\langle 4, g+1, g+3\rangle=\langle 4, 4n+3, 4n+5\rangle$. Since $4$, $4n+3$, and $4n+5$ are pairwise relatively prime, we apply Lemma~\ref{lem:three-mutually-coprime} to conclude that $S$ is not a free numerical semigroup.

Suppose $\mult(S)=4$ and $g\equiv3\pmod*{4}$. Then $g=4n+3$ for some $n\in\N_0$, and $S=\langle 4, g+2, g+3, 2g+1\rangle=\langle 4, 4n+5, 4n+6,8n+7\rangle$. However, $\gcd(4,8n+7)=\gcd(4n+5,4n+6)=1$. We apply Lemma~\ref{lem:distinct-coprime-pairs} to conclude that $S$ is not a free numerical semigroup.

If $\mult(S)=5$, then $g\equiv1,2,3,4\pmod*{5}$. In each case $A$, the minimal generating set for $S$, contains three of the four integers in $\ints{g+1,g+4}$. Thus, $A$ contains two consecutive integers $g+i$ and $g+i+1$ for some $i\in\ints{1,3}$. Note that $5\nmid g+i, g+i+1$. We also have $5\in A$. Then $\gcd(5,g+i)=\gcd(5,g+i+1)=\gcd(g+i,g+i+1)=1$, and so by Lemma~\ref{lem:three-mutually-coprime}, we conclude that $S$ is not a free numerical semigroup.

If $\mult(S)=6$, then $g\equiv1,2,3,4,5\pmod*{6}$. If $g\equiv3\pmod*{6}$, for $A$ the minimal generating set of $S$, we have $g+1,g+2,g+4,g+5\in A$. Then $\gcd(g+1,g+2)=\gcd(g+4,g+5)=1$, and so by Lemma~\ref{lem:distinct-coprime-pairs}, we conclude that $S$ is not a free numerical semigroup. If $g\not\equiv3\pmod*{6}$, then $A$ contains three consecutive integers which are congruent to either $1, 2, 3$ or $-3,-2,-1$ $\pmod*{6}$. In either case, these three integers are pairwise relatively prime. By Lemma~\ref{lem:three-mutually-coprime}, $S$ is not a free numerical semigroup.

If $\mult(S)\ge7$, let $a=\mult(S)$. We have $g\equiv1,2,\dots,a-1\pmod*{a}$. If $g\not\equiv\pm1,\pm2\pmod*{a}$, for $A$ the minimal generating set of $S$, we have  $g+1,g+2,g+(a-2),g+(a-1)\in A$. Since $\gcd(g+1,g+2)=\gcd(g+(a-2),g+(a-1))=1$, by Lemma~\ref{lem:distinct-coprime-pairs}, we conclude that $S$ is not a free numerical semigroup. If $g$ is 1 or 2 modulo $a$, then $g+(a-4),g+(a-3),g+(a-2),g+(a-1)\in A$. If $g$ is $-1$ or $-2$ modulo $a$, then $g+1,g+2,g+3,g+4\in A$. Thus, if $g\equiv \pm1,\pm2\pmod*{a}$, $A$ contains a subsequence of three consecutive integers which begins with an odd integer. Hence the greatest common divisor of each pair of these three integers is 1. By Lemma~\ref{lem:three-mutually-coprime}, $S$ is not a free numerical semigroup.
\end{proof}

\subsection{Wilf's Conjecture}
We will show that reflective numerical semigroups satisfy Wilf's conjecture, as stated in \cite{Wilf1979}.
\begin{cor}
    Suppose $S$ is a reflective numerical semigroup with $\genus(S)=g$ and $\mult(S)=a$. Then $\embdim(S)\ge\dfrac{\frob(S)+1}{\frob(S)+1-\genus(S)}$.
\end{cor}
\begin{proof}
Suppose $S$ is reflective with $\genus(S)=g$ and $\mult(S)=a$. Then $a\nmid g$. As usual, let $q,r\in\Z$ be such that $g=aq+r$ with $r\in\ints{1,a-1}$. By Corollary~\ref{cor:ME-frobenius}, $\frob(S)=2g-r$. We have 
\[
    \dfrac{\frob(S)+1}{\frob(S)+1-\genus(S)}
    = \dfrac{2g-r+1}{2g-r+1-g}
    = \dfrac{2qa+2r-r+1}{qa+r-r+1}
    = \dfrac{2qa+r+1}{qa+1}.
\]

The value of $\embdim(S)$ depends on $r$. 
If $r=a-1$, then $\embdim(S)=a$, and we have
\[
    \dfrac{2qa+r+1}{qa+1}
    =\dfrac{2qa+a}{qa+1}
    =2+\dfrac{a-2}{qa+1}
    \le 2+\dfrac{a-2}{1}
    =a
    =\embdim(S).
\]
If $r\in\ints{1,a-2}$, then $\embdim(S)=a-1$, and we have
\[
    \dfrac{2qa+r+1}{qa+1}
    = 2 + \dfrac{r-1}{qa+1}
    \le 2+\dfrac{r-1}{1}
    = r+1 
    \le a-1 
    = \embdim(S).
\]

Thus, if $S$ is reflective, then $\dfrac{\frob(S)+1}{\frob(S)+1-\genus(S)}\le\embdim(S)$.
\end{proof}

\section{Counting reflective numerical semigroups by genus and Frobenius number}\label{sec:counting-ME-semigroups}
Finally, we consider the problem of computing the number of reflective numerical semigroups that have a given genus or a given Frobenius number. Our motivation comes from the problem of computing the number of numerical semigroups of a given genus $g$ or of a given Frobenius number $F$. For recent results on these counts and related problems which remain open, see \cite{Bras-Amoros08} as well as \cite{Kaplan2017} and the references cited within.

As was mentioned in the introduction, we will make use of the multiplicative function $\tau(n)$, which, for $n\in\N$, is equal to the number of positive divisors of $n$.

The material in this section relies on the following corollary, which itself follows directly from Proposition~\ref{prop:ED iff cases}.

\begin{cor}\label{cor:at-most-one-reflective-semigroup}
    Let $g,a\in\N$. If $a\in\ints{2,g+1}$ and $a\nmid g$, then there is one reflective numerical semigroup $S$ with $\genus(S)=g$ and $\mult(S)=a$. If $a\not\in\ints{2,g+1}$ or $a\mid g$, then there are no reflective numerical semigroups $S$ with $\genus(S)=g$ and $\mult(S)=a$.
\end{cor}

\subsection{Counting by genus}
For $g\in\N$, let $\numg(g)$ equal the number of reflective numerical semigroups with genus $g$, and let $\numgsym(g)$ equal the number of such numerical semigroups which are also symmetric. 

In Example~\ref{ex:count-reflective-genus-6}, we found that there are three reflective numerical semigroups of genus 6. We can take the same approach to count the number of reflective numerical semigroups of genus $g$.  

\begin{prop}\label{prop:number-of-ME-semigroups-by-genus}
    Let $g\ge1$. Then $\numg(g) = g + 1 - \tau(g)$.
\end{prop}
\begin{proof}
    Fix $g\ge1$. By Corollary~\ref{cor:at-most-one-reflective-semigroup}, $\numg(g)=\#\{a\in\ints{2,g+1} : a\nmid g\}$. Since 1 is a divisor of $g+1$, We can replace the interval $\ints{2,g+1}$ with $\ints{1,g+1}$ without changing the result. Thus,
    \begin{align*}
        \numg(g) &= \#\{a\in\ints{2,g+1} : a\nmid g\} \\
        &= \#\{a\in\ints{1,g+1} : a\nmid g\} \\
        &= \#\ints{1,g+1} \setminus \{d\ge1 : d\mid g\} \\
        &= (g+1)-\tau(g),
    \end{align*}
    as desired.
\end{proof}

We now count the number of symmetric reflective numerical semigroups of genus $g$.

\begin{figure}
    \centering 
    \begin{tabular}{|c|c||c|c||c|c||c|c||c|c|} \hline 
    $g$ & $\mathrm{ng}(g)$ & $g$ & $\mathrm{ng}(g)$ & $g$ & $\mathrm{ng}(g)$ & $g$ & $\mathrm{ng}(g)$ & $g$ & $\mathrm{ng}(g)$ \\ \hline \hline 
    1 & 1 & 11 & 10 & 21 & 18 & 31 & 30 & 41 & 40  \\ \hline
    2 & 1 & 12 & 7 & 22 & 19 & 32 & 27 & 42 & 35  \\ \hline
    3 & 2 & 13 & 12 & 23 & 22 & 33 & 30 & 43 & 42  \\ \hline
    4 & 2 & 14 & 11 & 24 & 17 & 34 & 31 & 44 & 39  \\ \hline
    5 & 4 & 15 & 12 & 25 & 23 & 35 & 32 & 45 & 40  \\ \hline
    6 & 3 & 16 & 12 & 26 & 23 & 36 & 28 & 46 & 43  \\ \hline
    7 & 6 & 17 & 16 & 27 & 24 & 37 & 36 & 47 & 46  \\ \hline
    8 & 5 & 18 & 13 & 28 & 23 & 38 & 35 & 48 & 39  \\ \hline
    9 & 7 & 19 & 18 & 29 & 28 & 39 & 36 & 49 & 47  \\ \hline
    10 & 7 & 20 & 15 & 30 & 23 & 40 & 33 & 50 & 45  \\ \hline
    \end{tabular}
    \caption{The number of reflective numerical semigroups of genus $g$ for $g\in\ints{1,50}$}
    \label{fig:genus-table}
\end{figure}

\begin{prop}\label{prop:num-symm-MEs}
    Let $g\ge1$. Then
\[\numgsym(g) = \begin{dcases*}
1 & if $g=1$,\\ 
\tau(g-1)-1 & if $g>1$.
\end{dcases*}\]
\end{prop}
\begin{proof}
By Corollary~\ref{cor:at-most-one-reflective-semigroup}, for each combination of $g,a\in\N$ with $a\in\ints{2,g+1}$ and $a\nmid g$, there is one reflective numerical semigroup $S$ for which $\genus(S)=g$ and $\mult(S)=a$. There are no other reflective numerical semigroups. By Corollary \ref{cor:ME and symm or pseudo-symm}, a reflective numerical semigroup is symmetric if and only if $g\equiv1\pmod*{a}$.

Thus, for a given genus $g$, counting the number of reflective numerical semigroups which are symmetric amounts to counting the number of integers $a$ for which $a\in\ints{2,g+1}$, $a\nmid g$, and $a\mid(g-1)$. 

If $g=1$, then we must have $a=2$. Thus, $\numgsym(1)=1$. If $g>1$, then the positive divisors of $g-1$ are all contained in $\ints{1,g+1}$. The only divisor of $g-1$ which also divides $g$ is 1, so the number of suitable values of $a$ is $\numgsym(g)=\tau(g-1)-1$.
\end{proof}

One such reflective and symmetric numerical semigroup is $S=\langle 3,11\rangle$, visualized in Figure~\ref{fig:semigroup-rotate-1}. Note that $\mult(S)=3$,  $g=\genus(S)=(3-1)(11-1)/2=10$, and $10\equiv1\pmod*{3}$.

In particular, when $g-1$ is prime, there is a unique symmetric reflective numerical semigroup of genus $g$. We describe each such semigroup below.

\begin{cor}
For $p$ prime, there is exactly one symmetric reflective numerical semigroup with $\genus(S)=p+1$. This semigroup is $S=\langle A\rangle$ for
\[
    A=\begin{dcases*}
        \{2,7\} & if $p=2$; and \\
        \{p\}\cup\ints{p+2,2p-1} & if $p$ is odd.
    \end{dcases*}
\]
\end{cor}
\begin{proof}
    Let $g=p+1$. By Proposition \ref{prop:num-symm-MEs}, since $\tau(g-1)-1=\tau(p)-1=1$, there is one reflective symmetric numerical semigroup $S$ with $\genus(S)=g=p+1$. Let $a=\mult(S)$. Since $a\mid (g-1)=p$ and $a>1$, we must have $a=p$. We use Proposition \ref{prop:min gen sets} to find a minimal generating set for $S$. 
    
    For $a=p=2$, we have $g=3\equiv-1\pmod*{a}$. Then $S$ is minimally generated by 
    \[
        A=\{a,2g+1\}=\{2,7\}.
    \]
    In this case, the gaps of $S=\langle A\rangle$ are $\{1,3,5\}=\{1,g,2g-1\}$, which is equidistributed modulo $g=3$.
    
    For $a=p>2$, we have $g=p+1\equiv1\pmod*{a}$ and $g\not\equiv-1\pmod*{a}$. Then $S$ is minimally generated by 
    \[
        A =\{a\}\cup\ints{g+1,g+(a-1)}= \{p\}\cup\ints{p+2,2p-1}.
    \]
    In this case, the gaps of $S=\langle A\rangle$ are $\ints{1,p-1}\cup\{p+1, 2p+1\}=\ints{1,g-2}\cup\{g,2g-1\}$, which is equidistributed modulo $g$.
\end{proof}

\subsection{Counting by Frobenius number}
For $F\in\N$, let $\numF(F)$ denote the number of reflective numerical semigroups with Frobenius number $F$. In this subsection, we will find exact and asymptotic formulas for $\numF(F)$.

\begin{lem}\label{lem:at-most-one-g-for-each-m}
    Given $F,a\in\N$, there is at most one reflective numerical semigroup $S$ with $\frob(S)=F$ and $\mult(S)=a$.
\end{lem}
\begin{proof}
    By Corollary~\ref{cor:at-most-one-reflective-semigroup}, for any combination of $g,a\in\N$, we know there is at most one reflective numerical semigroup $S$ for which $\genus(S)=g$ and $\mult(S)=a$. It will therefore suffice to show that, for any combination of $a,F\in\N$, a reflective numerical semigroup $S$ with $\frob(S)=F$ and $\mult(S)=a$, there is a unique $g\in\N$ for which $\genus(S)=g$.
    
    By Corollary~\ref{cor:ME-frobenius}, if $S$ is reflective then $\frob(S)=2\genus(S)-r$, where $r$ is the unique integer in $\ints{1,\mult(S)-1}$ for which $\genus(S)\equiv r\pmod*{\mult(S)}$. Observe that 
    \[\frob(S)=2\genus(S)-r\equiv 2r-r\equiv r\pmod*{\mult(S)}\]
    as well. Thus, if we have $\frob(S)=F$ and $\mult(S)=a$, then we can compute $r\in\ints{1,a-1}$ such that $F\equiv r\pmod*{a}$. It follows that $\genus(S)=(F+r)/2$. Since this is the only possible value for $\genus(S)$, we are done.
\end{proof}    

\begin{prop}
    For $F,a\in\N$, there is a reflective numerical semigroup $S$ with $\frob(S)=F$ and $\mult(S)=a$ if and only if the following three conditions hold:
    \begin{enumerate} 
    \item $a\in\ints{2,F+1}$;
    \item $a\nmid F$; and
    \item $\Floor{F/a}\equiv0\pmod*{2}$.
    \end{enumerate}
\end{prop}
\begin{proof}
    As usual, we apply the division algorithm to $g$ and $a$ to obtain integers $q$ and $r$ for which $g=qa+r$ and $r\in\ints{1,a-1}$. (Recall that by Remark~\ref{rmk:g in s iff g=0}, $a\nmid g$.) 
    
    $(\implies)$. Suppose $S$ is reflective with $\frob(S)=F\ge1$, $\mult(S)=a$, and $\genus(S)=g$. Since $F\ge1$, we must have $a\ge2$. Since $F+1\in S$, we have $a\le F+1$. Next, since $a\in S$ and $F\not\in S$, $F$ cannot be a multiple of $a$. This means $a\nmid F$. We have shown the first two conditions.
    
    Next, by Corollary~\ref{cor:ME-frobenius}, $F=2g-r$. Then $F=2qa+r$. Dividing both sides by $a$ and taking the floor of each side, we see that $\Floor{F/a}=\Floor{2q+r/a}=2q$, an even integer.
    
    $(\impliedby)$. Let $a,F\in\N$ be such that $a\in\ints{2,F+1}$, $a\nmid F$, and $\Floor{F/a}=2k$ for some integer $k$. We will split up the first condition into three cases: $a=F+1$; $F/2\le a\le F$; or $2\le a<F/2$. We'll consider these cases separately.
    
    If $a=F+1$, consider the numerical semigroup $S=\langle\ints[a,2a-1]\rangle$, which has $\mult(S)=a$, $\HH(S)=\ints{1,a-1}$, $\genus(S)=a-1$, and $\frob(S)=a-1$. Note that $\HH(S)$ has one gap in each congruence class modulo $\genus(S)$, which means $S$ is reflective. Further, $\frob(S)=a-1=F$. Thus, there is a reflective numerical semigroup $S$ with $\mult(S)=a=F+1$ and $\frob(S)=F$.
    
    If $F/2\le a\le F$, then $\Floor{F/a}=1$. This violates the second condition on $a$, so we will not consider such a value of $a$.
    
    If $a\in\ints{2,F/2}$, then $a\le F/2$. Since $\Floor{F/a}$ is even and $a\nmid F$, $F=2ka+r$ for some integers $k,r$ with $r\in\ints{1,a-1}$. Let $g=(F+r)/2$. Then $g=ka+r$, an integer, and we have $a\le F/2 = g-r/2 <g$. Therefore $a\in\ints{2,g}$. And since $g=ka+r$ with $r\in\ints{1,a-1}$, $a\nmid g$. By Corollary~\ref{cor:at-most-one-reflective-semigroup}, since $a\in\ints{2,g}\subset\ints{2,g+1}$ and $a\nmid g$, there is a reflective numerical semigroup of genus $\genus(S)=g$ and multiplicity $\mult(S)=a$. By Corollary~\ref{cor:ME-frobenius}, this reflective numerical semigroup has Frobenius number $\frob(S)=2\genus(S)-r=2g-r=F$.
\end{proof}

We are now able to count the number of reflective numerical semigroups with Frobenius number $F$ by counting the number of integers $a$ in $\ints{2,F+1}$ that do not divide $F$ and for which $\Floor{F/a}$ is even. We note two things. First, there is always a reflective numerical semigroup with $a=F+1$, so we can count that as 1 in our formula for $\numF(F)$ and just check values of $a$ in $\ints{2,F}$. Second, since 1 divides $F$ and we only consider $a$ for which $a\nmid F$, we change nothing by considering $a\in\ints{1,F}$ instead of $\ints{2,F}$.
\begin{cor}\label{cor:nFrob-formula}
    For $F\in\N$, the number of reflective numerical semigroups $S$ with $\frob(S)=F$ is given by 
    \[
        \numF(F) 
        = 1 + \#\left\{a\in\ints{1,F} : a\nmid F,\, \Floor{F/a}\equiv0\pmod*{2}\right\}.
    \]
\end{cor}

Next, for $n\in\N$, let $\EE(n)=\{a\in\ints{1,n} : \Floor{n/a}\equiv0\pmod*{2}\}$. Additionally, following from the definition of $\tau(n)$ as the number of positive divisors of $n$, we let $\tau_e(n)$ denote the number of even positive divisors of $n$.

\begin{lem}\label{lem:n(F)-and-E(F)}
    For $F\in\N$, 
    $\numF(F)=1+\#\EE(F)-\tau_e(F)$.
\end{lem}
\begin{proof}
    By Corollary~\ref{cor:nFrob-formula}
    \[
        \numF(F)=1+\#\left\{a\in\ints{1,F} : a\nmid F, \Floor{F/a}\equiv0\pmod*{2}\right\}.
    \]
    To evaluate this, we can count each integer $a\in\ints{1,F}$ for which $\Floor{F/a}$ is even and then remove each such $a$ for which $a\mid F$. The number of integers $a$ for which $\Floor{F/a}$ is even is $\#\EE(F)$. The number of such integers $a$ for which $a\mid F$ and $F/a$ is even is exactly the number of even divisors of $F$. Therefore, 
    $\numF(F)=1+\#\EE(F)-\tau_e(F)$.
\end{proof}

To evaluate $\#\EE(F)$, we need to count the number of even terms in the sequence 
\begin{equation}
    \label{seq:fraction-floors}
    \Floor{F/1}, \Floor{F/2}, \dots, \Floor{F/F}.
\end{equation}
In \cite{DentShor2022}, a formula and an asymptotic formula are given for the number of odd terms in Sequence~\eqref{seq:fraction-floors}. The asymptotic formula, which follows from Dirichlet's solution to the \emph{Dirichlet divisor problem}, is given using big O notation. (See, e.g., \cite[Chapter 3]{Apostol1976}.)

\begin{prop}[{\cite[Proposition 2.2, Proposition 2.6]{DentShor2022}}]\label{prop:dentshor}
    For $n\in\N$,
    \[
        \#\left\{k : 1\le k\le n,\, \Floor{n/k}\equiv1\pmod*{2} \right\} 
        = \sum\limits_{k=1}^n (-1)^{k+1}\Floor{n/k} = n\log2 + \bigO\left(\sqrt{n}\right),
    \]
    where $\log$ denotes a logarithm using base $\ee$.
\end{prop}
This gives the number of odd terms in Sequence~\eqref{seq:fraction-floors}. If we subtract the above quantity from $n$, then we have the number of even terms in Sequence~\eqref{seq:fraction-floors}. Combining Lemma~\ref{lem:n(F)-and-E(F)} and Prop.~\ref{prop:dentshor}, we get a formula for $\numF(F)$.
\begin{cor}\label{cor:nF-exact}
    For $F\in\N$, $\numF(F)=1-\tau_e(F)+\sum\limits_{k=2}^F(-1)^{k}\Floor{F/k}$.
\end{cor}

\begin{figure}
    \centering
    \includegraphics[width=.45\textwidth]{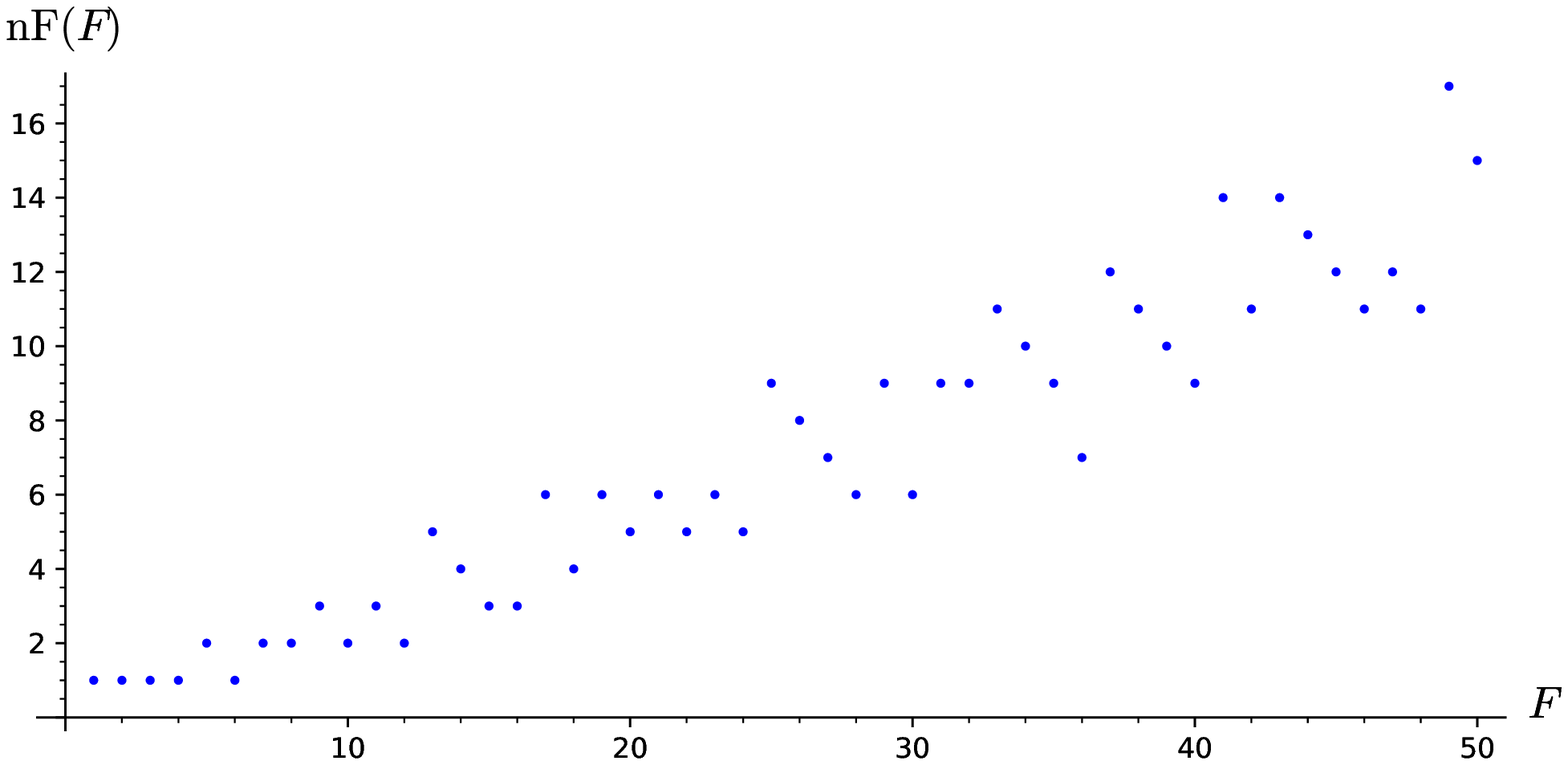}
    \includegraphics[width=.45\textwidth]{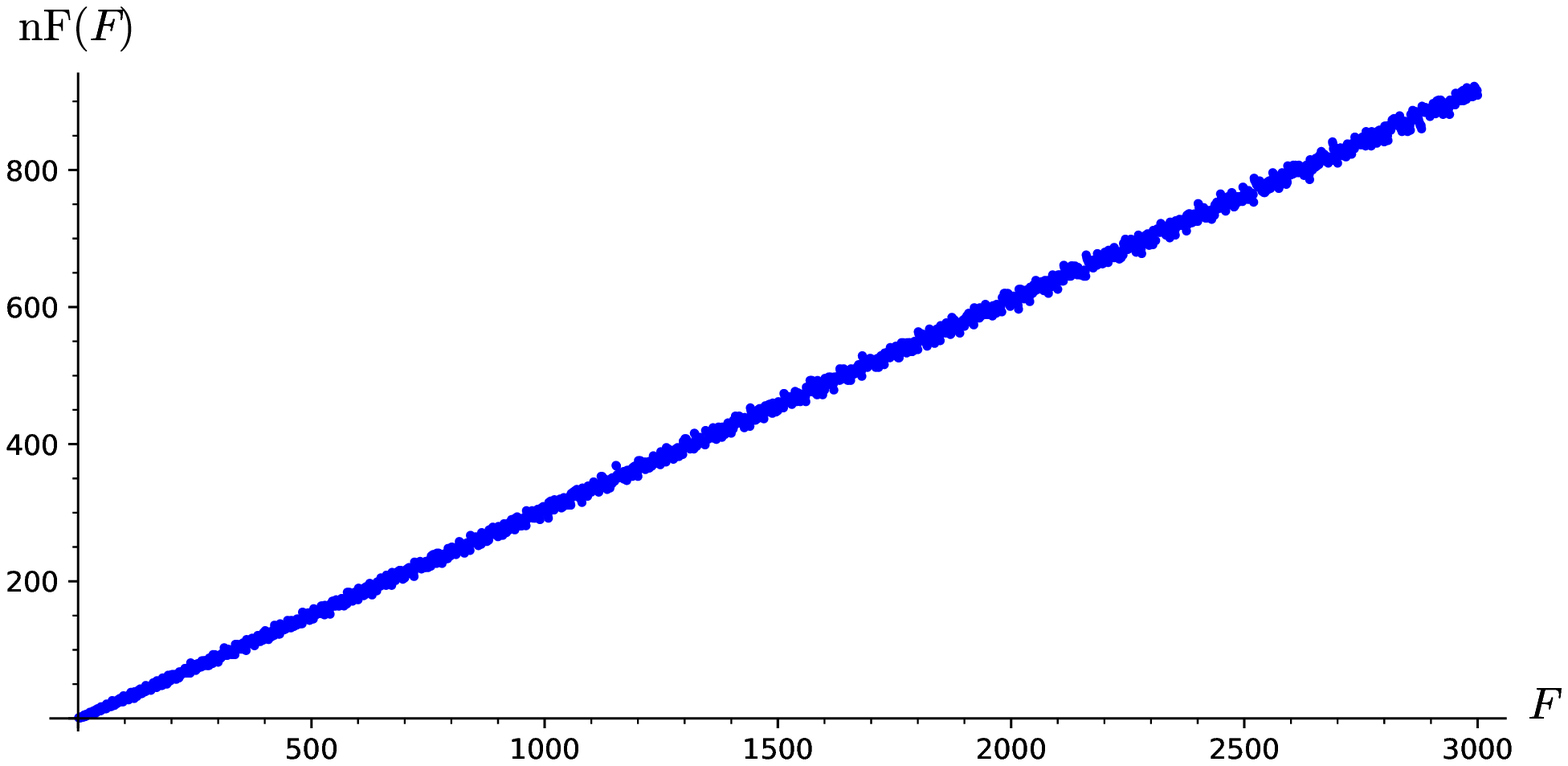}
    \caption{Plots of $\numF(F)$, the number of reflective numerical semigroups $S$ with $\frob(S)=F$, for $F\in\ints{1,50}$ and for $F\in\ints{1,3000}$.}
    
    \label{fig:reflective_frob_plot}
\end{figure}

Values of $\numF(F)$, for $1\le F\le 50$, can be found in Figure~\ref{fig:frobenius-table}. Plots of $\numF(F)$, for $1\le F\le 50$ and for $1\le F\le 3000$, which were created with Sage \cite{Sage}, appear in Figure~\ref{fig:reflective_frob_plot}.

\begin{figure}
    \centering
    \begin{tabular}{|c|c||c|c||c|c||c|c||c|c|} \hline 
$F$ & $\mathrm{nF}(F)$ & $F$ & $\mathrm{nF}(F)$ & $F$ & $\mathrm{nF}(F)$ & $F$ & $\mathrm{nF}(F)$ & $F$ & $\mathrm{nF}(F)$\\ \hline \hline 
1 & 1 & 11 & 3 & 21 & 6 & 31 & 9 & 41 & 14  \\ \hline
2 & 1 & 12 & 2 & 22 & 5 & 32 & 9 & 42 & 11  \\ \hline
3 & 1 & 13 & 5 & 23 & 6 & 33 & 11 & 43 & 14  \\ \hline
4 & 1 & 14 & 4 & 24 & 5 & 34 & 10 & 44 & 13  \\ \hline
5 & 2 & 15 & 3 & 25 & 9 & 35 & 9 & 45 & 12  \\ \hline
6 & 1 & 16 & 3 & 26 & 8 & 36 & 7 & 46 & 11  \\ \hline
7 & 2 & 17 & 6 & 27 & 7 & 37 & 12 & 47 & 12  \\ \hline
8 & 2 & 18 & 4 & 28 & 6 & 38 & 11 & 48 & 11  \\ \hline
9 & 3 & 19 & 6 & 29 & 9 & 39 & 10 & 49 & 17  \\ \hline
10 & 2 & 20 & 5 & 30 & 6 & 40 & 9 & 50 & 15  \\ \hline
\end{tabular}
    \caption{The number of reflective numerical semigroups with Frobenius number $F$ for $F\in\ints{1,50}$}
    \label{fig:frobenius-table}
\end{figure}

Finally, we can obtain an asymptotic formula for $\numF(F)$. All we need is an asymptotic formula for $\tau_e(F)$.
\begin{lem}
    For $n\in\N$, $\tau_e(n)=\bigO\left(\sqrt{n}\right)$.
\end{lem}
\begin{proof}
    Let $n\in\N$. If $n=bc$ for $b,c\in\N$ with $b\le c$, then $1\le b\le \sqrt{n}$ and $\sqrt{n}\le c\le n$. There are at most $\sqrt{n}$ such divisors $b$. Hence $\tau(n)\le 2\sqrt{n}$. Since $0\le \tau_e(n)<\tau(n)$, we have $\tau_e(n)\le2\sqrt{n}$, implying $\tau_e(F)=\bigO\left(\sqrt{n}\right)$.
\end{proof}

We conclude with an asymptotic formula for $\numF(F)$.
\begin{prop}\label{prop:nF-asymp}
    For $F\in\N$, $\numF(F)=(1-\log2)F+\bigO(\sqrt{F})$ and, therefore,
    \[\lim\limits_{F\to\infty}\frac{\numF}{F}=1-\log2\approx0.306853.\]
\end{prop}

Thus, for large $F$, $\numF(F)\approx(1-\log2)F$. We can see this roughly linear behavior in Figure~\ref{fig:reflective_frob_plot}.

\bibliography{refs} 
\bibliographystyle{abbrvnat}
\end{document}